\documentclass[final]{siamltex}
\usepackage{amsmath,epsf}
\usepackage{amssymb}
\usepackage{graphicx}
\usepackage{epsfig,epsf,latexsym,subfigure,cite}
\usepackage{color}
\usepackage{siunitx}
\newtheorem{algorithm}{{\bf  Scheme}}%[section]

\renewcommand\theequation{\thesection.\arabic{equation}}

\newcommand{\bm}{\boldsymbol}
\newcommand{\bu}{{\bm u}}
\newcommand{\bv}{{\bm v}}

\newcommand{\assign}{:=}

\newcommand{\q}{{U}}

\def\BX{\bm{X}}
\def\BY{\bm{Y}}

\newcommand{\bx}{{\bm x}}

\newcommand{\Bu}{{\bm u}}
\newcommand{\nhalf}{{{n+\frac 12}}}
\newcommand{\none}{{{n+1}}}

\newcommand{\wtilde}{\widetilde}

\newcommand{\Bn}{{\bm n}}

\newcommand{\Grad}[1]{\nabla #1}

\newcommand{\veps}{\varepsilon}
\newcommand{\dd}{\mathrm{d}}

\newtheorem{remark}{\textbf{Remark}}[section]

\newcommand{\be}{\begin{equation}}
\newcommand{\ee}{\end{equation}}

\newcommand{\bse}{\begin{subequations}}
\newcommand{\ese}{\end{subequations}}

\begin{document}

\title{Efficient Second Order {Unconditionally} Stable
Schemes for a Phase-field Moving Contact Line Model 
{Using Invariant Energy Quadratization Approach}}

\author{Xiaofeng Yang\thanks{
		Department of Mathematics, University of South Carolina,
		Columbia, SC, USA 29208. ({\tt xfyang@math.sc.edu}); and School of Mathematical Sciences, University of Electronic Science and Technology of China, Chengdu 610054, China.
		}
	\and
	Haijun Yu\thanks{Corresponding author.
		NCMIS \& LSEC, Institute of Computational Mathematics and Scientific/Engineering Computing,
		Academy of Mathematics and Systems Science, Beijing 100190, China;
		School of Mathematical Sciences, University of Chinese Academy of Sciences, Beijing  100049, China.
		(\tt hyu@lsec.cc.ac.cn)}	}

\maketitle

\begin{abstract}
  We consider the numerical approximations for a phase field
  model consisting of incompressible Navier-Stokes equations
  with a generalized Navier boundary condition, and the
  Cahn-Hilliard equation with a dynamic moving contact line
  boundary condition. A crucial and challenging issue for
  solving this model numerically is the time marching
  problem, due to the high order, nonlinear and coupled
  properties of the system. We solve this issue by
  developing two linear, second-order accurate and energy stable
  schemes based on the projection method for the
  Navier-Stokes equations, the invariant energy
  quadratization for the nonlinear gradient terms in the bulk and
  boundary, and a subtle implicit-explicit treatment for the
  stress and convective terms. The well-posedness of the
  semi-discretized system and the unconditional energy
  stabilities are proved. Various numerical results based on
  a spectral-Galerkin spatial discretization are presented
  to verify the accuracy and efficiency of the proposed
  schemes.
\end{abstract}

\begin{keywords} 
  Cahn-Hilliard equation, moving contact line,
  unconditionally energy stable, high order scheme,
  invariant energy quadratization
\end{keywords}

\begin{AMS}
	65M12, 65Z05, 65P40
\end{AMS}

\section{Introduction} 
In this paper, we consider numerical approximations for a
hydrodynamically-coupled phase field model
\cite{QianWS2003,QianWS2006} with moving contact
line (MCL) boundary conditions. Phase field method is a
popular approach that is widely used to simulate the
interfacial dynamics of multiple material components due to
its versatility in modeling as well as numerical simulations.
The fluid-fluid interface in this
method is considered as a continuous but steep change of
some physical property of two fluid components, e.g.,
density or viscosity, etc. An order parameter (phase field variable) is introduced to label the two fluid
components, the interface is then represented by a thin
but smooth transition layer that can remove the
singularities in practice. The standard phase field model
for incompressible immiscible fluid mixture is a nonlinear
system that couples the Cahn-Hilliard equation and the
Navier-Stokes equations via convective and stress terms.
Once the fluid-fluid interface touches a solid wall, a MCL
problem is induced. This phenomenon exists in many physical
and engineering processes such as wetting, coating, or
painting, etc. In this situation, the no-slip boundary
condition for the Navier-Stokes equations is no longer
applicable (cf. \cite{Moffatt1964,
  DussanD1974,Dussan1979}). Simulations in
\cite{KoplikBW1988,KoplikBW1989,ThompsonR1989} using
molecular dynamics simulations showed that nearly complete
slip happens near the MCL. In the context of phase field
method, a set of accurate boundary conditions for the MCL
problem was derived by Qian et.\,al.\,in
\cite{QianWS2003,QianWS2006}, resulting in a
standard two phase model supplemented with a generalized
Navier boundary condition (GNBC) and a dynamic contact line
condition (DCLC), where the nonlinear couplings also show up
in the boundary conditions.

Numerically, although the phase field variable is continuous
and smooth, the model is still very stiff since a small
dimensionless parameter related to the thickness of the
interface layer, is involved, for which certain numerical
methods like {fully implicit or explicit methods for
  nonlinear terms, need a very small time marching step size
  ({see the analysis} in
  \cite{Feng03,Xu06,ShenYang2010})}. Hence a challenging
issue for solving the model is to develop numerically stable
methods, namely, to establish efficient numerical schemes
that can verify the so-called energy stable property at the
discrete level, irrespective of the coarseness of the
discretization. Such a kind of algorithms is usually called
unconditionally energy stable or thermodynamically
consistent. Schemes with this property are specially
preferred since it is critical for the numerical schemes to
use larger time steps to capture the correct long time
dynamics that can reduce the time cost in
computing. Nonetheless, we need mention a basic fact that a
larger time step will definitely induce larger numerical
errors. In other words, schemes with unconditional energy
stability can allow arbitrarily large time step only for the
sake of the stability concern. To measure whether a scheme
is reliable or not, the controllable accuracy is another
important factor in addition to the stability. Therefore, if
one attempts to use a time step as large as possible while
maintaining the desirable accuracy, the only possible choice
is to develop more accurate schemes, e.g., second order
energy stable schemes, that is the main focus of this paper.

It is remarkable that, unlike the enormous numerical scheme
developments on the standard phase field model for two phase
fluid flows system with easy boundary conditions (without
MCLs), e.g., see
\cite{Liu.S03,FLSY05, Kim05, Mie10,SuLiXu10, Hua11, Li&W2013}, 
almost all
developed time marching schemes for solving the model with
MCLs are first order, e.g., see
\cite{HeGW2011,DongS2012,Dong2012, GaoW2012,
  Salgado2013, AlandC2014, GaoW2014, Shen.YY2015, Macicp2017,YuYang2017}. 
More precisely, to the
best of the authors' knowledge, no schemes can be claimed to
posses the following three properties, namely, {\it
  linear}, {\it unconditionally energy stable}
and {\it second order accuracy}. This is because two
additional numerical difficulties, the discretization of the
GNBC and DCLC conditions on the boundary, emerge besides the
regular stiffness issue induced by the nonlinear double well
potential in the Cahn-Hilliard equation.  At the very least,
even for the Cahn-Hilliard equation, the algorithm design is
still challenging.  
To overcome the stiffness issue, many efforts had
been implemented to remove the time step constraint,
including, e.g., the nonlinear convex splitting approach
\cite{Eyre98,Wan11,SWW11,Kim2014}, and the linear
stabilization approach
\cite{YFLS06, ShenYang2010, Yang14,ShenYSINUM2015,
  Shen.YY2015, CJYZ15, LSY15, Noc16, SY_density}. 
About the pros and cons
of these two methods, we give some detailed discussions in
Section 3.

Therefore, the aim of this paper is to develop some more
efficient and accurate schemes for solving the phase-field
{MCL} model proposed in \cite {QianWS2003, QianWS2006}. We
shall construct second order time stepping schemes which
satisfy a discrete energy law by combining several
successful approaches including the {\em projection method}
for the Navier-Stokes equations to decouple the velocity and
pressure, the {\em invariant energy quadratization} (IEQ)
method (cf. \cite{YANG02,YANG01,YANG03,YANG04,YZW2017}) for
nonlinear gradient terms that appear in the bulk as well as the
boundary of the phase field equation, and a subtle {\em
  implicit-explicit} treatment for the stress and convective
terms. At each time step, one can solve a linear elliptic
system for the phase variable and the velocity field, and a
Poisson equation for the pressure. We shall give a rigorous
proof of the well-posedness of the linear system together
with numerical results to verify the second order accuracy
in time and the efficiency.

The rest of the paper is organized as follows. In Section 2,
we briefly describe the phase field model with MCL boundary
conditions and its associated energy dissipation law.  In
Section 3, we present the numerical schemes, and prove the
well-posedness of the semi-discretized linear system and
their discrete energy dissipation law rigorously. In Section
4, we describe the implementation based on a
Fourier-Legendre spectral-Galerkin spatial
discretization. In section 5, we present various numerical
examples to illustrate the accuracy and efficiency of the
proposed schemes. Some concluding remarks are given in
Section 6.

\section{The phase field MCL model and its energy law}

Consider a mixture of two
immiscible, incompressible fluids in a confined domain
$\Omega\subset \mathbb {R}^d$ ($d=2,3$) with matched density
and viscosity. We introduce a phase field variable
 $\phi(\bx, t)$ such that
\begin{align}\label{hhhh}
\phi(\bx,t)=\begin{cases} 1, & \text{fluid I},\\
-1, & \text{fluid II}.
\end{cases}
\end{align}
%To have a thin, smooth transition region of width $O(\veps)$,
We consider the following Ginzburg-Landau free energy for the mixture
\begin{align}
E_{mix}(\phi)= \lambda\int_\Omega 
\Big(\frac {\veps}{2}|\nabla \phi|^2+F(\phi)\Big) d\bx,
\end{align}
where $\lambda$ denotes rescaled characteristic strength of
phase mixing energy. The first gradient term in $E_{mix}$
contributes to the hydro-philic type (tendency of mixing) of
interactions between the materials and the second part, the
double well bulk energy $F(\phi)=\frac
1{4\veps}(\phi^2-1)^2$, represents the hydro-phobic type
(tendency of separation) of interactions. As the consequence
of the competition between the two types of interactions,
the equilibrium configuration will include a diffusive
interface with a thickness proportional to the parameter
$\veps$.

The total bulk energy of the hydrodynamic system is a sum of the
kinetic energy $E_k$ together with the mixing energy
$E_{mix}$:
\begin{eqnarray}
E_{bulk}(\bu,\phi) = E_k(\bu) + E_{mix}(\phi) =\int_\Omega 
\left(\frac{1}{2}|\bu|^2+\lambda\left(\frac{\veps}{2}|\Grad\phi|^2+F(\phi) \right)\right)d\bx.
\end{eqnarray}
Here $\bu$ is the fluid velocity field, and we assume the fluid density is $1$.

The evolution of the phase function is governed by the
Cahn-Hilliard equation
\begin{eqnarray}
&& \phi_t+\Grad\cdot(\bu\phi)= M \Delta \mu,\label{phaseCH:1}\\
&&\mu= \lambda\left(-\veps\Delta \phi+f(\phi)\right),\label{phaseCH:2}
\end{eqnarray}
where $\mu$ is the chemical potential, $M$ is a
mobility parameter, and
$f(\phi)=F'(\phi)=\frac{1}{\veps}\phi(\phi^2-1)$. 
The momentum equation for the
hydrodynamics takes the usual form of the Navier-Stokes
equation as follows (cf. e.g. \cite{Liu.S03,Kim05, QianWS2006, Shen.YY2015})
%  GaoW2012,SYW12,LSY15,Li&W2013,Boyer06,Boy10}),
\begin{eqnarray}
&&\bu_t+(\bu\cdot\Grad)\bu-\nu\Delta\bu+\Grad p+\phi \Grad\mu=0,\label{eq:NS:1}\\
&&\Grad\cdot\bu=0,\label{eq:NS:2}
\end{eqnarray}
where $p$ is the pressure, $\nu$ is the kinetic viscosity of
the mixture. 

On the boundary $\Gamma$, we use the GNBC
for the velocity as follows \cite{QianWS2003,QianWS2006}, 
\begin{eqnarray} 
&&\Bu \cdot \Bn =0,\label{eq:nsbc}\\
&&\nu\partial_{\Bn}\Bu_\tau=-\nu \ell(\phi) (\Bu_\tau - \Bu_w)  
-\frac{\lambda}{\gamma} \dot{\phi}\nabla_{\tau} \phi, \label{eq:GNBC}
\end{eqnarray} 
and together with the DCLC for the phase field variable on the boundary $\Gamma$,
\begin{eqnarray} 
&&\partial_\Bn \mu\ = 0,  \label{eq:chbcpot}\\
&&\veps\partial_\Bn\phi=-\frac{1}{\gamma}\dot{\phi}-g(\phi),\label{eq:DCLC}
\end{eqnarray}
where $\dot{\phi}=\phi_t+\bu_\tau\cdot\Grad_\tau\phi$,
$\ell(\phi)\geq 0$ is a given coefficient function that is
the ratio of domain length to the slip length, $\gamma$ is a
boundary relaxation coefficient, $\bu_w$ is the wall
velocity, $\bu_\tau$ is the tangential velocity along the
boundary tangential direction $\tau$,
$\nabla_{\tau}=\nabla-(\Bn\cdot\nabla)\Bn$ is the gradient
along $\tau$, $g(\phi)=G'(\phi)$ and $G(\phi)$ is the
interfacial energy that is defined as 
\begin{equation} \label{eq:cle}
G(\phi) = - \frac{\sqrt{2}}{3} \cos\theta_s \sin\left( \frac{\pi}{2} \phi \right),
\end{equation}
where $\theta_s$ is the static contact angle. From
(\ref{eq:nsbc}), we have $\Bu = \Bu_{\tau}$ on boundary
$\Gamma$.

We now describe the energy dissipation law for the above
model. Here and after, for any function $f,g \in
L^{2}({\Omega})$, we use $(f,g)=\int_{\Omega}f(\bx) g(\bx)
\dd \bx$ to denote the $L^2$ inner product between functions
$f(\bx)$ and $g(\bx)$, $(f,g)_{\Gamma}$ to denote
$\int_{\Gamma}f(s) g(s) \dd s$, and $\|f\|^{2}=(f,f)$ and
$\|f\|^{2}_{\Gamma}=(f,f)_{\Gamma}$.

\begin{lemma} \label{thm:PDE} The  Navier-Stokes-Cahn-Hilliard (NSCH) system with GNBC and DCLC
(\eqref{phaseCH:1}-\eqref{eq:NS:2} with \eqref{eq:nsbc}-\eqref{eq:DCLC}) satisfies the following energy
dissipation law
\begin{equation} \label{eq:nschelaw}
\frac{d}{d t} E(\Bu,\phi) = - \nu\| \Grad\Bu \|^2 
- M \| \nabla \mu \|^2 
- \frac{\lambda }{\gamma} \| \dot{\phi} \|^2_{\Gamma} 
-\nu \| \sqrt{\ell(\phi)} \Bu_s \|^2_{\Gamma} 
-\nu\left(\ell(\phi) \Bu_s, \Bu_w \right)_{\Gamma},
\end{equation} 
where $\Bu_s = \Bu - \Bu_w$ is the velocity slip on
boundary $\Gamma$, and
\begin{equation} \label{eq:energies}
E(\Bu,\phi) = \int_\Omega\Big(\frac 12 |\Bu|^2+ \lambda\big(\frac{\veps}{2}| \nabla \phi |^2 
+F(\phi)\big)\Big)d\bx+\lambda\int_\Gamma G(\phi(s))ds.
\end{equation}
\end{lemma}

\begin{proof}
	The theorem is identical to Theorem 1 in \cite{Shen.YY2015} if 
	we let $L(\phi)=-\dot{\phi}/\gamma $
\end{proof}

\section{Numerical schemes}

We now construct time marching schemes to solve the NSCH
system
\eqref{phaseCH:1}-\eqref{phaseCH:2}-\eqref{eq:NS:1}-\eqref{eq:NS:2}
with boundary conditions of GNBC
\eqref{eq:nsbc}-\eqref{eq:GNBC} and DCLC
\eqref{eq:chbcpot}-\eqref{eq:DCLC}. With the aim of
constructing schemes that are linear, second order, and
unconditionally energy stable, we notice that there are
several numerical challenges, including (i) how to decouple
the computations of velocity and pressure; (ii) how to
discretize $f(\phi)$; (iii) how to discretize $g(\phi)$; and
(iv) how to develop proper discretizations for convective
and stress terms.

The first difficulty (i) actually has been well studied
during the last forty years, e.g., the projection type
methods are one of the best ways to solve it (cf.\,the
review in \cite{GMS06} and the references therein). The
difficulty (ii) is also well studied recently by two class
of methods where one is the nonlinear convex splitting
method (cf. e.g. \cite{Eyre98,Wan11,SWW11,Kim2014}), and the other is
the linear stabilization approach(cf. e.g. 
\cite{YFLS06, ShenYang2010, Yang14,ShenYSINUM2015,
	Shen.YY2015, CJYZ15, LSY15, Noc16, SY_density}).
  The convex splitting approach is energy stable, however,
  it produces nonlinear schemes at most cases, thus the
  implementations are often complicated and the
  computational costs are high. The linear stabilization
  approach introduces purely linear schemes, thus it is easy
  to implement.  But, its stability usually requests a
  special property (generalized maximum
  principle)\cite{ShenYang2010,WangYu2017} satisfied by the
  classical PDE solution and the numerical solution, that is
  very hard to prove in general. {Recently, some theoretical
    progress has been made to overcome this barrier for first order
    \cite{li_characterizing_2016} and second
    order~\cite{li_second_2017} stabilization methods, using subtle Fourier analysis.} However,
  to prove the unconditional stability, large
  stabilization constants are required.
In this paper, we use a newly developed IEQ
approach, which has been successfully applied to solve
several gradient flow type models
(cf. \cite{YANG04,YANG02,YANG01,YANG03,YZWS2017,YZW2017,YL2017}). Its
idea is to make the free energy quadratic in terms of new
variables via the change of variables. Then the free energy
and the PDE system are transformed into equivalent forms and
thus the nonlinear terms can be treated semi-explicitly.

More precisely, we define two new variables
\begin{eqnarray}\label{new:variable}
U=\phi^2-1, \quad W=\sqrt{G(\phi)+C},
\end{eqnarray}
where $C$ is a constant to ensure $G(\phi)+C$ be positive,
for instance, $C=\frac{\sqrt{2}}{3}+\eta$ with any
$\eta>0$. Hence we can rewrite the bulk and total free energy as
\begin{eqnarray} \label{new:def:eneb}
&&E_{bulk}(\Bu,\phi, U) = \int_\Omega\Big(\frac 12 |\Bu|^2+ \lambda\big(\frac{\veps}{2}| \nabla \phi |^2 
+\frac{1}{4\veps}U^2\big)\Big)d\bx,\\
%\end{equation}\vskip-0.3cm
&& \label{new:def:ene}
E(\Bu,\phi, U, W) = E_{bulk}(\Bu,\phi, U)
+\lambda\int_\Gamma W^2 ds-\lambda C|\Gamma|.
\end{eqnarray}
Thus, we have an equivalent new PDE system as follows
\begin{eqnarray}
&&\phi_t+\Grad\cdot(\bu\phi)=M\Delta \mu, \label{sys:1}\\
&&\mu=\lambda (-\veps \Delta\phi+\frac{1}{\veps}\phi U),\label{sys:2}\\
&&\bu_t+(\bu\cdot\Grad)\bu-\nu\Delta\bu+\Grad p+\phi\Grad\mu=0,\label{sys:3}\\
&&\Grad\cdot\bu=0,\label{sys:4}\\
&&U_t=2\phi\phi_t,\label{sys:5}
\end{eqnarray}
with the GNBC on $\Gamma$ as
\begin{eqnarray} 
&&\Bu \cdot \Bn= 0,\label{gn:1}\\
&& \nu\partial_{\Bn}\Bu_\tau =-\nu \ell(\phi) (\Bu_\tau - \Bu_w)
-\frac{ \lambda}{\gamma} \dot{\phi}\nabla_{\tau} \phi,\label{gn:2}
\end{eqnarray} 
and DCLC as
\begin{eqnarray}
&&\partial_\Bn \mu\ = 0,  \label{dcl:1}\\
&& \veps\partial_\Bn\phi=-\frac{1}{\gamma}\dot{\phi}-Z(\phi)W,\label{dcl:2}\\
&&W_t=\frac 12 Z(\phi)\phi_t,  \label{dcl:3}
\end{eqnarray}
where $Z(\phi)={g(\phi)}/{\sqrt{G(\phi)+C}}$.
The initial conditions read as
\begin{eqnarray}
\phi|_{t=0}=\phi_0,\quad 
\bu|_{t=0}=\bu_0,\quad
U|_{t=0}=\phi_0^2-1,\quad
W|_{t=0}=\sqrt{G(\phi^0)+C}.
\end{eqnarray}

We derive the energy dissipative law for this new system
\eqref{sys:1}-\eqref{dcl:3} as follows.

\begin{theorem} \label{thm2:PDE} The NSCH system with GNBC and
  DCLC \eqref{sys:1}-\eqref{dcl:3} satisfies the
  following energy dissipation law
\begin{eqnarray} \label{new:law}
\begin{aligned}
\frac{d}{d t} E(\bu,\phi, U,W) =& - \nu\| \Grad\Bu \|^2 
- M \| \nabla \mu \|^2 
- \frac{\lambda}{\gamma} \|\dot{\phi} \|^2_{\Gamma}
-\nu \| \sqrt{\ell(\phi)} \Bu_s \|^2_{\Gamma} 
\\
&
-\nu\left(\ell(\phi) \Bu_s, \Bu_w \right)_{\Gamma},
\end{aligned}
\end{eqnarray} 
%where $E_{tot}(\bu,\phi, U,W)$ is defined in \eqref{new:def:ene}.
\end{theorem}
\begin{proof}
  By taking the $L^2$ inner product of equation
  (\ref{sys:1}) with $\mu$, and using boundary conditions
  \eqref{gn:1} and \eqref{dcl:1}, we get
\begin{equation} \label{lem2:1}
\left( \phi_t, \mu \right) - \left( \Bu \phi, \nabla \mu \right) 
= - M\| \nabla \mu \|^2.
\end{equation}
By taking the $L^2$ inner product of equation (\ref{sys:2}) with
$- \phi_t $, we have
\begin{equation} \label{lem2:2}
- \left( \mu, \phi_t \right) 
= \lambda\left(\veps  \partial_\Bn \phi, \phi_t \right)_\Gamma
- \frac{\lambda\veps}{2} \frac{d}{dt} \| \nabla \phi \|^2 
- \frac{\lambda}{\veps} \left( \phi U, \phi_t\right).
\end{equation}
By taking the $L^2$ inner product of \eqref{sys:5} with
$\frac{\lambda}{2\veps}U$, we obtain
\begin{eqnarray}\label{lem2:3}
\lambda\frac{d}{dt}\frac{1}{4\veps}\|U\|^2=\frac{\lambda}{\veps}(\phi\phi_t, U).
\end{eqnarray}
By taking the $L^2$ inner product of equation (\ref{sys:3})
with $\Bu$, using the divergence free condition
(\ref{sys:4}), we get
\begin{eqnarray} \label{lem2:4}
\begin{aligned}
\frac{d}{d t} \frac 12\| \Bu \|^2=-\nu\|\Grad\bu\|^2 
+(\nu\partial_\Bn \Bu_\tau,\Bu_\tau)_\Gamma-(\phi\Grad\mu,\bu).
\end{aligned}
\end{eqnarray}
By taking the summation of \eqref{lem2:1}-\eqref{lem2:4}, we obtain
\begin{equation} \label{lem2:5}
\begin{aligned}
\frac{d}{dt} E(\bu,\phi, U)
= - \nu\| \Grad\Bu \|^2 - M\| \nabla \mu \|^2 
+ \left ( \nu \partial_\Bn \Bu, \Bu \right)_\Gamma + \lambda\left( \veps \partial_\Bn\phi, \phi_t \right)_\Gamma.
\end{aligned}
\end{equation} 
Then, we use boundary condition (\ref{gn:2}),
(\ref{dcl:2}) and definition of $\dot{\phi}$, to derive
\begin{eqnarray} \label{lem2:6}
\left( \nu\partial_{\Bn}\Bu, \Bu \right)_\Gamma 
&=& -\frac{\lambda}{\gamma} \Big( \dot{\phi} \nabla_{\tau} \phi, \Bu_{\tau} \Big)_\Gamma 
-\nu  \left(\ell(\phi) \Bu_s, \Bu_s + \Bu_w \right)_\Gamma,\\
\label{lem2:7}
\lambda (  \veps\partial_\Bn \phi, \phi_t)_\Gamma
&=&-\frac{\lambda}{\gamma}\|\dot{\phi}\|_\Gamma^2+\frac{\lambda}{\gamma}(\dot{\phi},\bu_\tau\cdot\Grad_\tau\phi)_\Gamma-\lambda(Z(\phi)W,\phi_t)_\Gamma.
\end{eqnarray}
By taking the $L^2$ inner product of \eqref{dcl:3} with $2\lambda  W$, we obtain
\begin{eqnarray}\label{lem2:8}
\lambda\frac{d}{dt}\|W\|_\Gamma^2=\lambda(Z(\phi)\phi_t,W)_\Gamma.
\end{eqnarray}

Summing up \eqref{lem2:5}-\eqref{lem2:8}, we get the desired energy law
(\ref{new:law}).
\end{proof}

  We emphasize that the new transformed system
  \eqref{sys:1}-\eqref{dcl:3} is exactly equivalent to the
  original system
  \eqref{phaseCH:1}-\eqref{phaseCH:2}-\eqref{eq:NS:1}-\eqref{eq:NS:2},
  \eqref{eq:nsbc}-\eqref{eq:GNBC}-\eqref{eq:chbcpot}-\eqref{eq:DCLC}
  that can be easily obtained by integrating \eqref{sys:5}
  and \eqref{dcl:3} with respect to the time. Therefore, the
  energy law \eqref{new:law} for the transformed system is
  exactly the same as the energy law \eqref{eq:energies} for
  the original system for the time-continuous case. We will
  develop time-marching schemes for the new transformed
  system \eqref{sys:1}-\eqref{dcl:3} that in turn follows
  the new energy dissipation law \eqref{new:law} instead of
  the energy law \eqref{eq:nschelaw} for the original
  system.

We fix some notations here. 
We define two Sobolev space $H_{c}^1(\Omega)=\{\phi\in
H^1(\Omega): \int_\Omega \phi=0\}$ and
$H_{\bu}(\Omega)=\big\{\bu\in [H^1(\Omega)]^d : \bu\cdot\Bn|_\Gamma=0\big\}$.
Let $\delta t >0$ be a time step size and set $t^n=n \delta
t$. For any function $S(\bx, t)$, 
let $S^n$ denotes the
numerical approximation to $S(\cdot,t)|_{t=t^n}$, and
$ S^{n+\frac 12 }:=\frac{S^{n+1}+S^n}{2}$, 
$S_\star^{n+\frac12}:=\frac32 S^{n}-\frac12 S^{n-1}$, 
$S_\star^{n+1}:=2S^{n}-S^{n-1}$, 
$\delta S^{n+1} :=S^{n+1}-S^n$, 
$\delta^2 S^{n+1} :=S^{n+1}-2S^n+S^{n-1}$. 

\subsection{A Crank-Nicolson scheme}

We first construct a linear Crank-Nicolson scheme (CN) for
the system \eqref{sys:1}-\eqref{dcl:3}, as follows.
\begin{algorithm}
  Assuming that $\phi^{n}$, $\bu^{n}$, $p^{n}$, $U^n$, $W^n$, $\phi^{n-1}$, $\bu^{n-1}$
  are given, we compute $\phi^{n+1}$, $\bu^{n+1}$,
  $p^{n+1}$, $U^{n+1}$, $W^{n+1}$ in two steps.

  {\bf Step 1:} We update $\phi^{n+1}, \mu^\nhalf,
  \wtilde\bu^{n+1}, U^{n+1}, W^{n+1}$ as follows,

\begin{eqnarray}
&&\frac{\phi^{n+1}-\phi^n}{\delta t}+
\Grad\cdot(\wtilde\bu^\nhalf\phi_\star^{n+\frac12})=M\Delta\mu^\nhalf,\label{cn:1}\\
&&\mu^\nhalf=\lambda\Big(-\veps\Delta\phi^\nhalf+\frac{1}{\veps}\phi_\star^{n+\frac12} U^\nhalf\Big),\label{cn:2}\\
&&U^{n+1}-U^n=2\phi_\star^{n+\frac12}(\phi^{n+1}-\phi^n),\label{cn:3}\\
&&\frac{\wtilde\bu^{n+1}-\bu^n}{\delta t}+B(\bu_\star^{n+\frac12},\wtilde\bu^\nhalf)-\nu\Delta\wtilde \bu^\nhalf+\Grad p^n+\phi_\star^{n+\frac12}\Grad\mu^\nhalf=0, \label{cn:4}
\end{eqnarray}
with the following boundary conditions on $\Gamma$,
\begin{eqnarray}
&&\wtilde\bu^{n+1}\cdot\Bn=0 ,\label{cnbd:1}\\
&&\nu\partial_{\Bn}\wtilde\bu_\tau^\nhalf=-\nu\ell(\phi_\star^{n+\frac12}) 
(\wtilde\bu^\nhalf -\bu_w)
-\frac{\lambda}{\gamma}\dot{\phi}^\nhalf \Grad_\tau\phi_\star^{n+\frac12} ,\label{cnbd:2}\\
&&\partial_\Bn\mu^\nhalf=0 ,\label{cnbd:3}\\
&&\veps\partial_\Bn\phi^\nhalf=-\frac{1}{\gamma}\dot{\phi}^\nhalf-Z(\phi_\star^{n+\frac12}) W^\nhalf,\label{cnbd:4}\\
&&W^{n+1}-W^n=\frac 12 Z(\phi_\star^{n+\frac12})(\phi^{n+1}-\phi^n), \label{cnbd:5}
\end{eqnarray}
where $B(\bu, \bv)=(\bu\cdot\Grad)\bv+\frac 12(\Grad\cdot\bu)\bv$,
$\wtilde\bu^\nhalf=\frac{\wtilde\bu^{n+1}+\bu^n}{2}$ and 
\begin{equation}
\dot{\phi}^\nhalf=\frac{\phi^{n+1}-\phi^n}{\delta t}+\wtilde\bu_\tau^\nhalf\cdot\Grad_\tau\phi_\star^{n+\frac12}.
\end{equation}

{\bf Step 2:} We update $\bu^{n+1}$ and $p^{n+1}$ as follows, 
\begin{eqnarray}
&&\frac{\bu^{n+1}-\wtilde\bu^{n+1}}{\delta t}+\Grad(\frac{p^{n+1}-p^n}{2})=0,\label{cn:5}\\
&&\Grad\cdot\bu^{n+1}=0,\label{cn:6}
\end{eqnarray}
with the boundary condition on $\Gamma$,
\begin{eqnarray}\label{cnbd:6}
\bu^{n+1}\cdot\Bn=0.
\end{eqnarray}

\end{algorithm}

\begin{remark}\label{rmk:pproj}
  The computations of $(\phi^{n+1},
  \mu^\nhalf,\wtilde\bu^{n+1})$ and the pressure $p^{n+1}$
  are totally decoupled via a second order pressure
  correction scheme \cite{vanKan1986} and a subtle
  implicit-explicit treatment for the stress and convective
  terms. It is quite an open problem on how to develop a
  second order scheme that can decouple the computations of
  $(\phi, \mu)$ from the velocity field $\bu$. All decoupled
  type energy stable schemes were first order accurate in
  time
  (cf. \cite{Yang14,LSY15,ShenYSINUM2015,Shen.YY2015,MinS13}). The
  adopted projection method here was analyzed in
  \cite{Shen1996} where it is shown (discrete time,
  continuous space) that the schemes is second order
  accurate for velocity in $\ell^2(0, T; L^2(\Omega))$ but
  only first order accurate for pressure in
  $\ell^\infty(0,T; L^2(\Omega))$. The loss of accuracy for
  pressure is due to the artificial boundary condition
  \eqref{cn:5} imposed on pressure \cite{ELi1995}.  We refer
  to \cite{Shen1996,GMS06} and references therein for
  analysis on this type of discretization.
\end{remark}

Schemes \eqref{cn:1}-\eqref{cnbd:5} is totally linear
since we handle the convective and stress terms by
compositions of implicit and explicit
discretization at $t^{n+\frac 12}$. 
Note that the new variables $U$ and $W$ will not bring up  extra computational cost. In fact, we can substitute for 
$U^{n+1}$ and $W^{n+1}$ in \eqref{cn:2} and \eqref{cnbd:4} 
using \eqref{cn:3} and \eqref{cnbd:5}.
Thus the scheme \eqref{cn:1}-\eqref{cnbd:5} can be written as 
\begin{eqnarray}\label{rewrite:cn1}
\left\{
\begin{aligned}
&\phi^{n+1}+\frac{\delta t}{2}\Grad\cdot(\wtilde\bu^{n+1}\phi_\star^{n+\frac12})-\delta tM\Delta\mu^\nhalf=f_1,\\
&-\mu^\nhalf-\frac{\lambda\veps}{2}\Delta\phi^\none+\frac{\lambda}{\veps}(\phi_\star^{n+\frac12})^2\phi^\none=f_2,\\
&\frac 12\wtilde\bu^{n+1}+\frac{\delta t}{4}B(\bu_\star^{n+\frac12},\wtilde\bu^{n+1})-\frac{\nu\delta t}{4}\Delta \bu^{n+1}+\frac{\delta t}{2}\phi_\star^{n+\frac12}\Grad\mu^\nhalf=f_3,
\end{aligned}
\right.
\end{eqnarray}
with the following boundary conditions on $\Gamma$, 
\begin{eqnarray}\label{rewrite:cn2}
\left\{
\begin{aligned}
&\wtilde\bu^{n+1}\cdot\Bn=0,\\
&\nu\partial_\Bn\wtilde\bu^{n+1}=-\nu\ell(\phi_\star^{n+\frac12})\wtilde\bu^{n+1}
-\frac{2\lambda}{\gamma}\mathring{\phi}^{ n+1}\Grad\phi_\star^{n+\frac12}+g_1,\\
&\partial_\Bn\mu^\nhalf=0,\\
&\veps\partial_\Bn\phi^{n+1}=-\frac{2}{\gamma}\mathring{\phi}^{n+1}-\frac 12 Z(\phi_\star^{n+\frac12})^2\phi^{n+1}+g_2,
\end{aligned}
\right.
\end{eqnarray}
where the definition of ${\mathring\phi}^{n+1}$ is 
\begin{eqnarray}
\hskip -2cm
\begin{aligned}
{\mathring\phi}^{n+1}=\frac{1}{\delta t}\phi^{n+1}+\frac{1}{2}\wtilde\bu^{n+1}_\tau\cdot\Grad_\tau\phi_\star^{n+\frac12},
\end{aligned}
\end{eqnarray}
and $f_1,f_2, f_3, g_1,g_2$ include only terms from previous time steps.
Therefore, we can solve
\eqref{rewrite:cn1}-\eqref{rewrite:cn2} directly. Once we
obtain $\phi^{n+1}, \mu^\nhalf, \wtilde \bu^{n+1}$, the new
variables $U^{n+1}$, $W^{n+1}$ are updated using
\eqref{cn:3} and \eqref{cnbd:5}.

Now we study the well-posedness of the semi-discretized system. 
Define $\bar\phi=\frac{1}{|\Omega|}\int_\Omega\phi d\bx$,
$\bar\mu=\frac{1}{|\Omega|}\int_\Omega\mu d\bx$. 
By integrating \eqref{cn:1}, we find that 
$\bar\phi^{n+1} = \bar\phi^{n} =\ldots = \bar\phi^0$.
Let $\phi=\phi^{n+1}-\bar\phi^0$, $\mu=\mu^\nhalf-\bar\mu^\nhalf$
such that $\phi, \mu \in H^1_c(\Omega)$. 
Here 
\begin{equation}\label{eq:muhalfb}
\bar\mu^\nhalf = \frac{1}{|\Omega|} \int_{\Omega}(-\frac{\lambda\veps}{2}\Delta\phi^\none+\frac{\lambda}{\veps}(\phi_\star^{n+\frac12})^2\phi^\none-f_2).
\end{equation}
Then the
weak form for \eqref{rewrite:cn1}-\eqref{rewrite:cn2} can be
written as the following system with unknowns $\mu,\phi \in H_{c}^1(\Omega)$, $\bu\in H_{\bu}(\Omega)$,
\begin{eqnarray}
&&(\phi,w)-\frac{\delta t}{2}(\bu\phi_\star^{n+\frac12},\Grad w)+M\delta t(\Grad \mu,\Grad w)=(f_1,w),\label{linear:1}\\
&&-(\mu,\psi)+\frac{\lambda\veps}{2}(\Grad\phi,\Grad\psi)+\frac{\lambda}{\veps}((\phi_\star^{n+\frac12})^2\phi,\psi)\nonumber\\
&&\hskip 0.5cm
+\frac{\lambda}{\gamma}(\mathring\phi,\psi)_\Gamma+\frac{\lambda}{4}(Z(\phi_\star^{n+\frac12})^2\phi,\psi)_\Gamma=(f_2,\psi)+\frac{\lambda}{2}(g_2,\psi)_\Gamma,\label{linear:2}\\
&&\frac{1}{2}(\bu,\bv)+\frac{\delta t}{4} (B(\bu_\star^{n+\frac12}, \bu),\bv)+\frac{\nu\delta t}{4}(\Grad\bu,\Grad\bv)+\frac{\delta t}{2}(\phi_\star^{n+\frac12}\Grad\mu,\bv)\nonumber\\
&&\hskip 0.5cm+\frac{\delta t\nu}{4}(\ell(\phi_\star^{n+\frac12})\bu,\bv)_\Gamma
+\frac{\lambda}{2\gamma}\delta t(\mathring\phi\Grad\phi_\star^{n+\frac12}, \bv)_\Gamma
=(f_3,\bv)
+\frac{\delta t}{4}(g_1,\bv)_\Gamma,\label{linear:3}
\end{eqnarray}
for any $w,\psi \in H_{c}^1(\Omega)$ and $\bv\in
H_{\bu}(\Omega)$, where $\mathring\phi=\frac{1}{\delta
  t}\phi+\frac 12\bu_\tau\cdot\Grad_\tau \phi_\star^{n+\frac12}$.

We denote the above linear system
\eqref{linear:1}-\eqref{linear:3} as
\begin{eqnarray}\label{final:matrix}
(\mathbb A \BX, \BY)= (\mathbb B, \BY),
\end{eqnarray}
where $\BX=(\mu,\phi,\bu)^T, \BY=(w,\psi,\bv)^T$ and $\BX, \BY\in (H_{c}^1,H_{c}^1, H_{\bu})(\Omega)$.

\begin{theorem} \label{wellpose} The linear system
  \eqref{linear:1}-\eqref{linear:3} admits a unique solution
  $(\mu, \phi, \bu)$ where $\mu, \phi\in
  H_{c}^1(\Omega)$, $\bu\in H_{\bu}(\Omega)$.
\end{theorem}
\begin{proof}
(i) For any $\BX=(\mu,\phi, \bu)^T$ and $\BY=(w,\psi, \bv)^T$ with $\BX, \BY\in(H_{c}^1,H_{c}^1, H_{\bu})(\Omega)$, we have
\begin{eqnarray}
(\mathbb A\BX, \BY)\le C_1(\|\phi\|_{H^1}+\|\mu\|_{H^1}+\|\bu\|_{H^1})(\|\psi\|_{H^1}+\|w\|_{H^1}+\|\bv\|_{H^1}),
\end{eqnarray}
from the trace theorem, where
$C_1$ is a constant depending on $\delta t$, $\nu$, $M$, $\veps$, $\lambda$, $\gamma$,
$\|\bu_\star^{n+\frac12}\|_\infty$, $\|\phi_\star^{n+\frac12}\|_\infty$,
$\|Z(\phi_\star^{n+\frac12})\|_\infty$. Therefore, the bilinear form
$(\mathbb A\BX, \BY)$ is bounded.

(ii) It is easy to derive that 
\begin{eqnarray}
\begin{aligned}
(\mathbb A\BX, \BX)=&\frac{1}{2}\|\bu\|^2+\frac{\nu\delta t}{4}\|\Grad \bu\|^2+\frac{\lambda\veps}{2}\|\Grad\phi\|^2+\frac{\lambda}{\veps}\|\phi_\star^{n+\frac12}\phi\|^2+M\delta t\|\Grad\mu\|^2\\
&+\frac{\lambda\delta t}{\gamma}\|\mathring\phi\|_\Gamma^2+\frac{\lambda}{4}\|Z(\phi_\star^{n+\frac12})\phi\|_\Gamma^2
+\frac{\delta t}{2}\nu\|(\ell(\phi_\star^{n+\frac12}))^{\frac 12}\bu_\tau\|_\Gamma^2\\
\ge& C_2(\|\phi\|_{H^1}^2+\|\mu\|_{H^1}^2+\|\bu\|_{H^1}^2), 
\end{aligned}
\end{eqnarray}
from Poincar\'e inequality
(since $\int_\Omega\phi d\bx=\int_\Omega\mu d\bx=0$), where $C_2$ is a constant depending on $\delta t, \nu, M,
\veps, \lambda$. Thus the bilinear form
$(\mathbb A\BX, \BY)$ is coercive.

Then from the Lax-Milgram theorem, we conclude the linear
system \eqref{final:matrix} admits a unique solution
$(\mu,\phi,\bu)\in (H_{c}^1, H_{c}^1,
H_{\bu})(\Omega)$. Namely, the linear system
\eqref{rewrite:cn1}-\eqref{rewrite:cn2} admits a unique
solution $(\mu^\nhalf, \phi^{n+1},\wtilde\bu^{n+1})$ in
$(H_{c}^1, H_{c}^1, H_{\bu})(\Omega)$.
\end{proof}

The stability result of the proposed Crank-Nicolson scheme
follows the same lines as in the derivation of the new PDE
energy dissipation law Theorem \ref{thm2:PDE}, as follows.
\begin{theorem}
The scheme \eqref{cn:1}-\eqref{cnbd:6} is unconditionally energy stable, in the sense that, it satisfies the following discrete
energy dissipation law,
\begin{eqnarray}\label{eq:ediss:cn2}
\begin{aligned}
E_{cn}^{n+1}=E_{cn}^n&-M\delta t\|\Grad\mu^\nhalf\|^2-\nu\delta t\|\Grad\wtilde\bu^\nhalf\|^2
-\frac{\lambda\delta t}{\gamma}\|\dot{\phi}^\nhalf\|_\Gamma^2,
\\
&
-\nu\delta t\|\ell(\phi_\star^{n+\frac12})^{\frac 12}\wtilde\bu_s^\nhalf\|_\Gamma^2
- \nu{\delta t} (\ell(\phi_\star^{n+\frac12})\wtilde\bu_s^\nhalf,\bu_w)_\Gamma
\end{aligned}
\end{eqnarray}
where $\wtilde\bu_s^\nhalf = \wtilde\bu^\nhalf - \bu_w$ and
\begin{eqnarray}\label{new:def:ene:disc}
\begin{aligned}
E_{cn}^n = E(\bu^n,\phi^n,U^n,W^n)+\frac{\delta t^2}{8}\|\Grad p^n\|^2.
\end{aligned}
\end{eqnarray}
%and $E_{tot}(\bu,\phi, U, W)$ is defined in \eqref{new:def:ene}.
\end{theorem}
\begin{proof}
  By taking the $L^2$ inner product of \eqref{cn:1} with
  $\delta t\mu^\nhalf$ and performing integration by parts,
  we obtain
\begin{eqnarray}\label{cnstab:1}
(\phi^{n+1}-\phi^n, \mu^\nhalf)-\delta t(\wtilde\bu^\nhalf\phi_\star^{n+\frac12},\Grad\mu^\nhalf)=-M\delta t\|\Grad\mu^\nhalf\|^2.
\end{eqnarray}
By taking the $L^2$ inner product of \eqref{cn:2} with $-(\phi^{n+1}-\phi^n)$, we obtain
\begin{eqnarray}\label{cnstab:2}
\begin{aligned}
-(\mu^\nhalf, \phi^{n+1}-\phi^n)=&
\lambda(\veps\partial_\Bn\phi^\nhalf,\phi^{n+1}-\phi^n)_\Gamma
-\frac{\lambda\veps}{2}(\|\Grad\phi^{n+1}\|^2-\|\Grad\phi^n\|^2)\\
&-\frac{\lambda}{\veps}(\phi_\star^{n+\frac12} U^\nhalf,\phi^{n+1}-\phi^n).
\end{aligned}
\end{eqnarray}
By taking the $L^2$ inner product of \eqref{cn:3} with $\frac{\lambda}{2\veps}U^\nhalf$, we obtain
\begin{eqnarray}\label{cnstab:3}
\begin{aligned}
\frac\lambda{4\veps}\big({\|U^{n+1}\|^2}-{\|U^n\|^2}\big)
=\frac{\lambda}{\veps}(\phi_\star^{n+\frac12}(\phi^{n+1}-\phi^n),U^\nhalf).
\end{aligned}
\end{eqnarray}
By taking the $L^2$ inner product of \eqref{cn:4} with
$\delta t\wtilde \bu^\nhalf$, we obtain
\begin{eqnarray}\label{cnstab:4}
\begin{aligned}
\frac 12 \|\wtilde\bu^{n+1}\|^2-\frac 12 \|\bu^n\|^2&+\delta t\nu\|\Grad\wtilde\bu^\nhalf\|^2-\delta t (\nu\partial_\Bn \wtilde\bu^\nhalf,\wtilde\bu^\nhalf)_\Gamma\\
&+\delta t(\Grad p^n,\wtilde\bu^\nhalf)+\delta t(\phi_\star^{n+\frac12}\Grad\mu^\nhalf,\wtilde\bu^\nhalf)=0.
\end{aligned}
\end{eqnarray}
By taking the $L^2$ inner product of \eqref{cn:5} with
$\delta t\bu^{n+1}$ and using the divergence free condition
for $\bu^{n+1}$ from \eqref{cn:6}, we obtain
\begin{eqnarray}\label{cnstab:5}
\begin{aligned}
\frac 12 (\|\bu^{n+1}\|^2-\|\wtilde\bu^{n+1}\|^2)+\frac 12\|\bu^{n+1}-\wtilde\bu^{n+1}\|^2=0.
\end{aligned}
\end{eqnarray}
We further rewrite the projection step \eqref{cn:5} as 
\begin{eqnarray}\label{cnstab:6}
\bu^{n+1}+\bu^n-2\wtilde\bu^\nhalf+\frac{\delta t}{2}\Grad(p^{n+1}-p^n)=0.
\end{eqnarray}
By taking the $L^2$ inner product of the above equation with
$\frac{\delta t}{2}\Grad p^n$ and applying the divergence
free condition for $\bu^{n+1}+\bu^n$, we obtain
\begin{eqnarray}\label{cnstab:7}
\frac{\delta t^2}{8}(\|\Grad p^{n+1}\|^2-\|\Grad p^n\|^2-\|\Grad p^{n+1}-\Grad p^n\|^2)=\delta t (\wtilde\bu^\nhalf,\Grad p^n).
\end{eqnarray}
On the other hand, it follows directly from \eqref{cn:5} that
\begin{eqnarray}\label{cnstab:8}
\frac{\delta t^2}{8}\|\Grad(p^{n+1}-p^n)\|^2=\frac{1}{2}\|\bu^{n+1}-\wtilde\bu^{n+1}\|^2.
\end{eqnarray}
Summing up \eqref{cnstab:1}, \eqref{cnstab:2}, \eqref{cnstab:3}, \eqref{cnstab:4}, \eqref{cnstab:5}, \eqref{cnstab:7} and \eqref{cnstab:8}, we obtain
\begin{eqnarray}\label{cnstab:9}
\begin{aligned}
&\frac{\lambda\veps}{2}\big(\|\Grad\phi^{n+1}\|^2- \|\Grad \phi^n\|^2\big)
+\frac\lambda{4\veps}\big({\|U^{n+1}\|^2}-{\|U^{n}\|^2}\big)
+M\delta t\|\Grad\mu^\nhalf\|^2\\
&\hskip1em +\frac 12 \|\bu^{n+1}\|^2-\frac 12\|\bu^n\|^2+\nu\delta t\|\Grad\wtilde\bu^\nhalf\|^2+\frac{\delta t^2}{8}(\|\Grad p^{n+1}\|^2-\|\Grad p^n\|^2)\\
&\hskip1em -\delta t(\nu\partial_\Bn\wtilde\bu^\nhalf,\wtilde\bu^\nhalf)_\Gamma-\lambda(\veps\partial_\Bn\phi^\nhalf,\phi^{n+1}-\phi^n)_\Gamma=0.
\end{aligned}
\end{eqnarray}
To deal with the boundary integrals,  from \eqref{cnbd:2}, we derive
\begin{eqnarray}\label{cnstab:10}
\begin{aligned}
-\delta t(\nu\partial_\Bn\wtilde\bu^\nhalf, \wtilde\bu^\nhalf)_\Gamma
=&\delta t\nu(\ell(\phi_\star^{n+\frac12})\wtilde\bu_s^\nhalf,\wtilde\bu_s^\nhalf+\bu_w)_\Gamma \\
%+ \delta t\nu(\ell(\phi_\star^{n+\frac12})\wtilde\bu_s^\nhalf,\bu_w)_\Gamma
&+\frac{\lambda\delta t}{\gamma}(\dot{\phi}^\nhalf\Grad_\tau\phi_\star^{n+\frac12},\wtilde\bu^\nhalf)_\Gamma.
\end{aligned}
\end{eqnarray}
From \eqref{cnbd:4} and the definition of $\dot{\phi}^\nhalf$, we obtain
\begin{eqnarray}\label{cnstab:11}
\begin{aligned}
&-\lambda(\veps\partial_\Bn\phi^\nhalf,\phi^{n+1}-\phi^n)_\Gamma
-\lambda(Z(\phi_\star^{n+\frac12}) W^\nhalf, \phi^{n+1}-\phi^n)_\Gamma\\
%&\hskip 1cm=\lambda(\frac{1}{\gamma}\dot{\phi}^\nhalf, %\phi^{n+1}-\phi^n)_\Gamma\\
&\hskip 2cm=\frac{\lambda\delta t}{\gamma}\|\dot{\phi}^\nhalf\|_\Gamma^2-\frac{\lambda\delta t}{\gamma}(\dot{\phi}^\nhalf,\wtilde\bu^\nhalf\cdot\Grad_\tau\phi_\star^{n+\frac12})_\Gamma.
\end{aligned}
\end{eqnarray}
By taking the $L^2$ inner product of \eqref{cnbd:5} with $2\lambda W^\nhalf$, we obtain
\begin{eqnarray}\label{cnstab:12}
\lambda(\|W^{n+1}\|_\Gamma^2-\|W^n\|_\Gamma^2)=\lambda (Z(\phi_\star^{n+\frac12}) (\phi^{n+1}-\phi^n), W^\nhalf)_\Gamma.
\end{eqnarray}
Finally, we obtain the desired energy law \eqref{eq:ediss:cn2} by
combining \eqref{cnstab:9}-\eqref{cnstab:12}.
\end{proof}

The proposed scheme follows the new energy dissipation law
\eqref{new:law} formally instead of the energy law for the
originated system \eqref{eq:nschelaw}. In the
time-continuous case, the two energy laws are the same. In
the time-discrete case, the discrete energy 
$E^n_{cn}$ (defined in
\eqref{new:def:ene:disc}) is a second order approximation
to the exact energy $E_{tot}(\bu^n,\phi^n)$ (defined in
\eqref{eq:energies}), since $U^{n+1},W^{n+1}$ are second
order approximations to $\phi^2-1$ and $\sqrt{G(\phi)+C} $.

Several remarks are in order.

\begin{remark}
  The time discretization for the cubic polynomial term
  $f(\phi)$ induced from the double well potential has been
  well-studied in a large quantity of literature.
  For instance, one popular method to obtain second
  order time marching schemes is the convex splitting
  approach (cf. \cite{Han2015,YANG01}) since there exists a
  natural convex-concave decomposition for the double well
  potential $F(\phi)$. For the boundary
  energy $G(\phi)$, since its second order derivative is
  bounded, it is natural to use the linear stabilization
  approach, see \cite{GaoW2014, GaoW2012,
    YuYang2017,Shen.YY2015}. Namely, $g(\phi)$ is treated
  explicitly and an extra linear stabilizer is added to
  improve the stability. However, it is not easy to design 
  unconditionally stable second order scheme by linear stabilization.

  The IEQ approach provides a novel way to handle both $f(\phi)$
  and $g(\phi)$. Its idea is very simple but quite different
  from the traditional time marching schemes like fully
  explicit, implicit or other various Taylor expansions to
  discretize nonlinear potentials. Through a simple
  substitution of new variables, the complicated nonlinear
  potentials are transformed into quadratic forms. We
  summarize the great advantages of this quadratic
  transformations as follows: (i) this quadratization method
  works well for various complex nonlinear terms as long as
  the corresponding nonlinear potentials are bounded from
  below; (ii) the complicated nonlinear potential is
  transferred to a quadratic polynomial form which is much
  easier to handle; (iii) the derivative of the quadratic
  polynomial is linear, which provides the fundamental
  support for linearization method; (iv) the quadratic
  formulation in terms of new variables can automatically
  maintain this property of positivity (or bounded from
  below) of the nonlinear potentials.
\end{remark}
\begin{remark}
  When the nonlinear potential is a fourth order
  polynomial, e.g., the double well potential, the IEQ we used in \eqref{new:variable} for $\phi$ variable  
   is exactly the same as the so-called Lagrange
  multiplier method developed in \cite{GuTi2013}. 
  We remark
  that the idea of Lagrange multiplier method only works
  well for the fourth order polynomial potential
  ($\phi^4$). This is because the nonlinear term $\phi^3$
  (the derivative of $\phi^4$) can be naturally decomposed
  into a multiplication of two factors: $\lambda(\phi)\phi$
  that is the Lagrange multiplier term, and the
  $\lambda(\phi)=\phi^2$ is then defined as the new
  auxiliary variable $U$. However, this method might not
  succeed for other type potentials. For instance, we notice
  the Flory-Huggins potential is widely used in two-phase
  model, see also \cite{Cahn&H1958}. The induced nonlinear
  term is logarithmic type as ${\rm ln}
  (\frac{\phi}{1-\phi})$. If one forcefully rewrites this
  term as $\lambda(\phi)\phi$, then $\lambda(\phi)=
  \frac{1}{\phi}{\rm ln} (\frac{\phi}{1-\phi})$ that is the
  definition of the new variable $U$. Obviously, such a form
  is unworkable for algorithms design. Therefore, we can see
  that the IEQ approach generalizes the Lagrange multiplier
  approach which is for double well potential only, and
  extends its applicability greatly to a unified framework
  for general dissipative stiff systems with high
  nonlinearity. About the application of the IEQ approach to
  handle other type of nonlinear potentials, we refer to the
  authors' other work in
  \cite{YANG02,YANG01,YANG03,YANG04,YZW2017}.

\end{remark}

\begin{remark} 
  The IEQ approach is more efficient than the nonlinear
  approach like fully implicit or convex splitting. Let us
  consider the double well potential case, e.g.,
  $F(\phi)=\phi^4$, then IEQ scheme will generate the linear
  scheme as $(\phi^{n})^2\phi^{n+1}$. The implicit or convex
  splitting approach will produces the scheme as
  $(\phi^{n+1})^3$. Therefore, if the Newton iterative
  method is applied for this term, at each iteration the
  nonlinear convex splitting approach would yield the same
  linear operator as IEQ approach. Hence the cost of solving
  the IEQ scheme is the same as the cost of performing one
  iteration of Newton method for the implicit/convex
  splitting approach, provided that the same linear solvers
  are applied.
\end{remark}

\begin{remark}
  {Instead of using $U=\phi^2-1$ as in equation
    \eqref{new:variable}, one can also use a more general
    form $U=\sqrt{(\phi^2-1)^2+C}$ with $C\ge0$. All the
    stability properties still hold. However the convergence
    behavior will be different. Our numerical tests show
    that the numerical results with $C\sim {\it o}(1)$
    perform almost the same as the results using
    \eqref{new:variable}. If $C$ is big, for a fixed time
    step, the magnitude of accuracy result is a little bit
    inferior to the case of $C=0$, but the accuracy order is
    still second order. Thus, for the double-well potential,
    one can either use $U=\phi^2-1$ or
    $U=\sqrt{(\phi^2-1)^2 + C}$ with $C\sim {\it
      o}(1)$.
    But, if the nonlinear potential is not double-well, for
    instance, the logarithmic Flory-Huggins potential, the
    only choice for the new variable $U$ is
    $\sqrt{G(\phi) + C}$ formally, see \cite{YANG04}.
    }
\end{remark}

\subsection{A backward differentiation scheme}
We further develop another linear scheme based
on second order backward differentiation formula(BDF2), that reads as follows.

\begin{algorithm}
  Assuming that $(\phi, \bu, p, U, W)^{n-1}$ and $(\phi,
  \bu, p, U, W)^{n}$ are already known, we compute
  $\phi^{n+1}, \bu^{n+1}, p^{n+1}, U^{n+1}, W^{n+1}$ from
  the following second order temporal semi-discrete system:

{\bf Step 1:} We update $\phi^{n+1}, \wtilde\bu^{n+1}, U^{n+1}, W^{n+1}$ as follows,
\begin{eqnarray}
&&\frac{3\phi^{n+1}-4\phi^n+\phi^{n-1}}{2\delta t}+
\Grad\cdot(\wtilde\bu^\none\phi_\star^{n+1})=
M\Delta\mu^\none,\label{bdf:1}\\
&&\mu^\none=\lambda\Big(-\veps\Delta\phi^\none+\frac{1}{\veps}\phi_\star^{n+1} U^\none\Big),\label{bdf:2}\\
&&3U^{n+1}-4U^n+U^{n-1}=2\phi_\star^{n+1}
(3\phi^{n+1}-4\phi^n+\phi^{n-1}),\label{bdf:3}\\
&&\frac{3\wtilde\bu^{n+1}\!\!\!-\!\!4\bu^n\!+\!\bu^{n\!-\!1}}{2\delta t}
+B(\bu_\star^{n+1}\!,\wtilde\bu^\none)\!-\!\nu\Delta\wtilde \bu^\none\!\!+\Grad p^n\!\!+\phi_\star^{n+1}\Grad\mu^\none=0, \label{bdf:4}
\end{eqnarray}
with the boundary conditions
\begin{eqnarray}
&&\wtilde\bu^{n+1}\cdot\Bn=0 ,\label{bdfbd:1}\\
&&\nu\partial_{\Bn}\wtilde\bu_\tau^\none=
-\nu\ell(\phi_\star^{n+1})(\wtilde\bu^\none -\bu_w)-\frac{\lambda}{\gamma}\dot{\phi}^\none \Grad_\tau\phi_\star^{n+1} ,\label{bdfbd:2}\\
&&\partial_\Bn\mu^\none=0 ,\label{bdfbd:3}\\
&&\veps\partial_\Bn\phi^\none=-\frac{1}{\gamma}\dot{\phi}^\none-Z(\phi_\star^{n+1}) W^\none,\label{bdfbd:4}\\
&&3W^{n+1}-4W^n+W^{n-1}=\frac 12 Z(\phi_\star^{n+1})(3\phi^{n+1}-4\phi^n+\phi^{n-1}), \label{bdfbd:5}
\end{eqnarray}
where 
\begin{eqnarray}
%\left\{
\begin{aligned}
\dot{\phi}^\none=\frac{3\phi^{n+1}-4\phi^n+\phi^{n-1}}{2\delta t}+\wtilde\bu_\tau^\none\cdot\Grad_\tau
\phi_\star^{n+1}.
\end{aligned}
%\right.
\end{eqnarray}

{\bf Step 2:} We update $\bu^{n+1}$ and $p^{n+1}$ as follows, 
\begin{eqnarray}
&&\frac3{2\delta t}\big(\bu^{n+1}-\wtilde\bu^{n+1}\big)
+\Grad(p^{n+1}-p^n)=0,\label{bdf:5}\\
&&\Grad\cdot\bu^{n+1}=0,\label{bdf:6}
\end{eqnarray}
with the boundary condition
\begin{eqnarray}\label{bdfbd:6}
\bu^{n+1}\cdot\Bn=0 \quad \mbox{on}\ \Gamma.
\end{eqnarray}

\end{algorithm}

Similar to the Crank-Nicolson scheme
\eqref{cn:1}-\eqref{cnbd:5}, the BDF2 scheme
\eqref{bdf:1}-\eqref{bdfbd:6} is linear and
the new variables $U, W$ do not
involve any extra computational cost.

\begin{theorem} \label{wellpose:bdf} The weak form of the linear
  system 
  \eqref{bdf:1}-\eqref{bdfbd:6}
  admits a
  unique solution $(\mu^\none$, $\phi^\none$,
  $\wtilde\bu^\none)$, where $\mu^\none,
  \phi^{n+1}\in H_{c}^1(\Omega)$, and $\wtilde\bu^{n+1}\in
  H_{\bu}(\Omega)$.
\end{theorem}
\begin{proof}
  The proof is similar to Theorem \ref{wellpose}, thus we
  omit the details here.
\end{proof}

\begin{theorem}
  The scheme \eqref{bdf:1}-\eqref{bdfbd:6} is
  unconditionally energy stable satisfying the following
  discrete energy dissipation law,
\begin{eqnarray}\label{ene:bdf2}
\begin{aligned}
E_{bdf}^{n+1}\le  E_{bdf}^n & 
-M\delta t\|\Grad\mu^{n+1}\|^2
-\nu\delta t\|\Grad\wtilde\bu^{n+1}\|^2
-\frac{\lambda}{\gamma}\delta t\|\dot{\phi}^{n+1}\|_\Gamma^2\\
&-\nu\delta t\|\ell(\phi_\star^{n+1})^{\frac 12}\wtilde\bu_s^{n+1}\|_\Gamma^2,
- \nu{\delta t} (\ell(\phi_\star^{n+1})\wtilde\bu_s^{n+1},\bu_w)_\Gamma,
\end{aligned}
\end{eqnarray}
where  
\begin{equation}\label{bdf:ene:def}
E_{bdf}^{n}= \frac12 E(\bu^n,\phi^n,U^n,W^n)
+ \frac12
E(\bu_\star^{n+1},\phi_\star^{n+1},U_\star^{n+1},W_\star^{n+1})
+\frac{\delta t^2}{3}\|\Grad p^n\|^2.
\end{equation}

\end{theorem} 
\begin{proof}
By taking the $L^2$ inner product of \eqref{bdf:1} with $2\delta t\mu^{n+1}$ and performing integration by parts, we obtain
\begin{equation}\label{stab2:0}
(3\phi^{n+1}\!\!-4\phi^n\!\!+\phi^{n-1}\!,\mu^{n+1})-2\delta t(\wtilde\bu^{n+1}\phi_\star^{n+1}\!,\Grad\mu^{n+1})=-2M\delta t\|\Grad\mu^\none\|^2.
\end{equation}
By taking the $L^2$ inner product of \eqref{bdf:2} with $-(3\phi^{n+1}-4\phi^n+\phi^{n-1})$, we obtain
\begin{eqnarray}\label{stab2:1}
\begin{aligned}
&-(\mu^\none, 3\phi^{n+1}-4\phi^n+\phi^{n-1})\\
&\hskip 0.5cm
=
-\frac{\lambda\veps}{2}\Big(\|\Grad\phi^{n+1}\|^2
-\|\Grad\phi^{n}\|^2
+\|\Grad\phi_\star^{n+2}\|^2
-\|\Grad\phi_\star^{n+1}\|^2
+\|\Grad\delta^2\phi^{n+1}\|^2\Big)\\
&\hskip 1cm +\lambda(\veps\partial_\Bn\phi^{n+1},3\phi^{n+1}-4\phi^n+\phi^{n-1})_\Gamma-\frac{\lambda}{\veps}(\phi_\star^{n+1} U^{n+1},3\phi^{n+1}-4\phi^n+\phi^{n-1}).
\end{aligned}
\end{eqnarray}
By taking the $L^2$ inner product of \eqref{bdf:3} with $\frac{\lambda}{2\veps}\q^{n+1}$, we obtain
\begin{equation}\label{stab2:2}
\hskip 1cm\begin{aligned}
\frac{\lambda}{4\veps}\Big(\|\q^{n+1}\|^2
&-\|\q^n\|^2
+\|\q_\star^{n+2}\|^2
-\|\q_\star^{n+1}\|^2
+\|\delta^2\q^{n+1}\|^2\Big)\\
&=
\frac{\lambda}{\veps}(\phi_\star^{n+1}(3\phi^{n+1}-4\phi^n+\phi^{n-1}), U^{n+1}).
\end{aligned}
\end{equation}
By taking the $L^2$ inner product of \eqref{bdf:4} with
$2\delta t\wtilde\bu^{n+1}$, we obtain
\begin{equation}\label{stab2:3}
\begin{aligned}
&(3\wtilde\bu^{n+1}\!\!-4\bu^n\!+\bu^{n-1}, \wtilde\bu^{n+1})+2\nu\delta t\|\Grad\wtilde\bu^{n+1}\|^2-2\nu\delta t(\partial_\Bn \wtilde\bu^{n+1},\wtilde\bu^{n+1})_\Gamma\\
&\hskip2.5cm +2\delta t(\Grad p^n,\wtilde\bu^{n+1})
+2\delta t(\phi_\star^{n+1}\Grad\mu^{n+1}, \wtilde\bu^{n+1})=0.
\end{aligned}
\end{equation}
From \eqref{bdf:5}, for any function $\bv$ with $\Grad\cdot \bv=0$, we can derive 
\begin{eqnarray}
(\bu^{n+1}, \bv)=(\wtilde\bu^{n+1}, \bv).
\end{eqnarray}
Then for the first term in \eqref{stab2:3}, we have
\begin{equation}\label{stab2:4}
\begin{aligned}
&(3\wtilde\bu^{n+1}\!\!-4\bu^n\!+\bu^{n-1}, \wtilde\bu^{n+1})\\
&\hskip 1cm=
(3\bu^{n+1}-4\bu^n+\bu^{n-1}, \bu^{n+1})
+3(\wtilde \bu^{n+1}-\bu^{n+1}, \wtilde \bu^{n+1}+\bu^{n+1})
\\
&\hskip 1cm=
\frac{1}{2}\Big(\|\bu^{n+1}\|^2-\|\bu^{n}\|^2
+\|\bu_\star^{n+2}\|^2
-\|\bu_\star^{n+1}\|^2
+\|\delta^2\bu^{n+1}\|^2\Big)\\
&\hskip1.4cm+ 3(\|\wtilde\bu^{n+1}\|^2-\|\bu^{n+1}\|^2).
\end{aligned}
\end{equation}

For the projection step, we rewrite \eqref{bdf:5} as
\begin{eqnarray}
\frac{3}{2\delta t}\bu^{n+1}+\Grad p^{n+1}=\frac{3}{2\delta t}\wtilde\bu^{n+1}+\Grad p^n.
\end{eqnarray}
By squaring both sides of the above equality, we obtain
\begin{equation}
\frac{9}{4\delta t^2}\|\bu^{n+1}\|^2+\|\Grad p^{n+1}\|^2=\frac{9}{4\delta t^2}\|\wtilde\bu^{n+1}\|^2+\|\Grad p^{n}\|^2+\frac{3}{\delta t}(\wtilde\bu^{n+1},\Grad p^n),
\end{equation}
namely, we have
\begin{equation}\label{stab2:5}
\frac{3}{2}(\|\bu^{n+1}\|^2-\|\wtilde\bu^{n+1}\|^2)+\frac{2\delta t^2}{3}(\|\Grad p^{n+1}\|^2-\|\Grad p^n\|^2)=2\delta t(\wtilde\bu^{n+1},\Grad p^n).
\end{equation}
By taking the $L^2$ inner product of \eqref{bdf:5} with $2\delta t\bu^{n+1}$, we have
 \begin{equation}\label{stab2:6}
\frac{3}{2}(\|\bu^{n+1}\|^2-\|\wtilde\bu^{n+1}\|^2+\|\bu^{n+1}-\wtilde\bu^{n+1}\|^2)=0.
\end{equation}
By combining \eqref{stab2:0}-\eqref{stab2:4} and \eqref{stab2:5}-\eqref{stab2:6}, we obtain
\begin{equation}\label{final:0}
\begin{aligned}
&2M\delta t\|\Grad\mu^{n+1}\|^2+\frac{3}{2}\|\bu^{n+1}-\wtilde\bu^{n+1}\|^2\\
&+\frac{2\delta t^2}{3}(\|\Grad p^{n+1}\|^2-\|\Grad p^n\|^2)+2\nu\delta t\|\Grad\wtilde\bu^{n+1}\|^2\\
&+\frac{\lambda\veps}{2}\Big(\|\Grad\phi^{n+1}\|^2
-\|\Grad\phi^{n}\|^2
+\|\Grad\phi_\star^{n+2}\|^2
-\|\Grad\phi_\star^{n+1}\|^2
+\|\Grad\delta^2\phi^{n+1}\|^2\Big)\\
&+\frac{\lambda}{4\veps}\Big(\|\q^{n+1}\|^2
-\|\q^{n}\|^2
+\|\q_\star^{n+2}\|^2
-\|\q_\star^{n+1}\|^2
+\|\delta^2\q^{n+1}\|^2\Big)\\
&+\frac{1}{2}\Big(\|\bu^{n+1}\|^2
-\|\bu^{n}\|^2
+\|\bu_\star^{n+2}\|^2
-\|\bu_\star^{n+1}\|^2
+\|\delta^2\bu^{n+1}\|^2\Big)\\
&=
\lambda(\veps\partial_\Bn\phi^{n+1},3\phi^{n+1}-4\phi^n+\phi^{n-1})_\Gamma
+2\delta t(\nu\partial_\Bn\wtilde\bu^{n+1},\wtilde\bu^{n+1})_\Gamma.
\end{aligned}
\end{equation}
From \eqref{bdfbd:4}, we obtain
\begin{equation}\label{stab2:7}
\begin{aligned}
&-\lambda\Big(\veps\partial_\Bn\phi^{n+1},
3\phi^{n+1}\!\!-4\phi^n\!+\phi^{n-1}\Big)_\Gamma
\\
&\hskip 2cm
=2\delta t\frac{\lambda}{\gamma}\|\dot{\phi}^{n+1}\|_\Gamma^2-2\delta t\frac{\lambda}{\gamma}(\dot{\phi}^{n+1},\wtilde\bu^{n+1}\cdot\Grad\phi_\star^{n+1})_\Gamma\\
&\hskip 2.5cm
+\lambda(Z(\phi_\star^{n+1}) W^{n+1},3\phi^{n+1}\!\!-4\phi^n\!+\phi^{n-1})_\Gamma.
\end{aligned}
\end{equation}
From \eqref{bdfbd:2}, we obtain
\begin{equation}\label{stab2:8}
\begin{aligned}
-\nu(\partial_\Bn\wtilde\bu^{n+1}, \wtilde\bu^{n+1})_\Gamma&
= \frac{\lambda}{\gamma}\big(\dot{\phi}^{n+1}\Grad\phi_\star^{n+1},\wtilde\bu^{n+1}\big)_\Gamma \\
&+\nu \|\ell(\phi_\star^{n+1})^\frac12\wtilde\bu_s^{n+1}\|^2_\Gamma
+\nu \big(\ell(\phi_\star^{n+1})\wtilde\bu_s^{n+1},\bu^w\big)_\Gamma.
\end{aligned}
\end{equation}
By taking the $L^2$ inner product of \eqref{bdfbd:5} with $2\lambda W^{n+1}$, we obtain
\begin{eqnarray}\label{stab2:9}
\begin{aligned}
&\lambda\Big(\|W^{n+1}\|^2_\Gamma
-\|W^n\|^2_\Gamma
+\|W_\star^{n+2}\|^2_\Gamma
-\|W_\star^{n+1}\|^2_\Gamma
+\|\delta^2W^{n+1}\|^2_\Gamma\Big)\\
&\hskip2cm=
\lambda \Big(Z(\phi_\star^{n+1}) (3W^{n+1}-4W^n+W^{n-1}), W^{n+1}\Big)_\Gamma.
\end{aligned}
\end{eqnarray}
By combining \eqref{final:0}, \eqref{stab2:7}-\eqref{stab2:9}, we obtain
\begin{equation}\label{eq:bdf:edissex}
\hskip-0.2cm\begin{aligned}
&2M\delta t\|\Grad\mu^{n+1}\|^2
+\frac{2\delta t^2}{3}(\|\Grad p^{n+1}\|^2-\|\Grad p^n\|^2)
+2\nu\delta t\|\Grad\wtilde\bu^{n+1}\|^2\\
&+\frac{\lambda\veps}{2}\Big(\|\Grad\phi^{n+1}\|^2-\|\Grad\phi^{n}\|^2
+\|\Grad\phi_\star^{n+2}\|^2
-\|\Grad\phi_\star^{n+1}\|^2
\Big)\\
&+\frac{\lambda}{4\veps}\Big(\|\q^{n+1}\|^2-\|\q^{n}\|^2
+\|\q_\star^{n+2}\|^2
-\|\q_\star^{n+1}\|^2
\Big)\\
&+\frac{1}{2}\Big(\|\bu^{n+1}\|^2-\|\bu^{n}\|^2
+\|\bu_\star^{n+2}\|^2
-\|\bu_\star^{n+1}\|^2
\Big)\\
&+\lambda\Big(\|W^{n+1}\|^2_\Gamma-\|W^n\|^2_\Gamma
+\|W_\star^{n+2}\|^2_\Gamma
-\|W_\star^{n+1}\|^2_\Gamma
\Big)\\
&+2\delta t\frac{\lambda}{\gamma}\|\dot{\phi}^{n+1}\|^2_\Gamma
+2\nu\delta t\|(\ell(\phi_\star^{n+1}))^\frac 12\wtilde\bu_s^{n+1}\|_\Gamma^2 
-2\nu\delta t (\ell(\phi_\star^{n+1})\wtilde\bu_s^\none, \bu_w)_\Gamma\\
&=
-\frac{3}{2}\|\bu^{n+1}\!-\wtilde\bu^{n+1}\|^2
\!-\|\Grad\delta^2\phi^{n+1}\|^2
\!-\|\delta^2\q^{n+1}\|^2
\!-\|\delta^2\bu^{n+1}\|^2
\!-\|\delta^2W^{n+1}\|^2_\Gamma\\
&\le 0.
\end{aligned}
\end{equation}
We conclude the theorem.
\end{proof}

\begin{remark} 
	{
		Same to the CN scheme, the BDF2 scheme is only first
		order accurate for pressure. As mentioned in Remark \ref{rmk:pproj}, the loss
		of accuracy for pressure is due to the artificial condition
		\begin{equation}
          \partial_\Bn p^{n+1} = 0, \quad \mbox{on}\ \Gamma. 
        \end{equation} imposed on pressure by \eqref{cn:5}-\eqref{cnbd:6} for CN scheme and 
        \eqref{bdf:5}-\eqref{bdfbd:6} for BDF2 scheme, respectively.
        The particularly semi-implicit treatment for the phase-field body force in Navier-Stokes equation allows 
        the pressure in a rotational form pressure-projection method satisfies an appropriate boundary condition. 
        For example, let's replace equation (3.71) with a rotational form (cf. \cite{timmermans_approximate_1996, guermond_error_2004, GMS06})
	\begin{equation*}
	\frac{3}{2\delta t}({\bu}^{n+1}-{\tilde{\bu}}^{n+1}) 
	+ \nabla (p^{n+1}-p^n + \nu \nabla\cdot \tilde{\bu})=0.
	\end{equation*}
    Noticing $ \phi_\star^{n+1}\partial_\Bn\mu^\none=0$ and
    $B(\bu_\star^{n+1}\!,\tilde\bu^\none)\cdot {\bf{n}}=0$
    on boundary, using similar procedure as in
    \cite{guermond_error_2004}, we find that $p$ satisfies
    boundary condition
	\begin{eqnarray*}
		\partial_\Bn p^{n+1} = -\nu (\nabla\times \nabla \times{\bu}^{n+1})\cdot {\bf{n}},\quad \mbox{on}\ \Gamma.	
	\end{eqnarray*}
    The body force due to phase field does not introduce any
    trouble because we use the form $\phi \nabla \mu$
    instead of $-\mu \nabla \phi$ which was frequently used
    in other literature.  Our numerical results show that
    the rotational form pressure-projection indeed improve
    the accuracy of pressure. But since our focus is not the
    accuracy of pressure, we will not discuss this issue
    in details.  }
\end{remark}

\section{Spatial discretization and implementation}

We implement a spectral Galerkin method for a 2-dimensional rectangular domain 
$\Omega= [0, L_x] \times [- 1, 1]$ to
test the stability and accuracy of the linear schemes proposed in last section.
Note that any spatial
discretization based on weak formulation and Galerkin
approximation of the NSCH coupled system will keep the
energy dissipation properties obtained for the temporal
semi-discretized schemes.  Thus the finite element or
spectral element methods can be used as well in a similar way.

Now we take the weak form for the Crank-Nicolson scheme
\eqref{linear:1}-\eqref{linear:3} as an example to
illustrate the spatial discretization, and the BDF2 scheme
can be handled similarly. We obtain a Galerkin approximation
by providing appropriate finite dimensional spaces
$H_{\Bu}^N, H_{\mu}^N, H_{\phi}^N$ for $H_{\Bu}(\Omega),
H^1_c(\Omega), H^1_c(\Omega)$ for the velocity field $\bu$,
the chemical potential $\mu$ and the phase variable $\phi$,
correspondingly.

 \subsection{A spectral Galerkin approximation}

 We assume the system in $x$ direction is periodic, only the
 top and bottom boundaries ($y=\pm 1$) take the GNBC and
 DCLC.  We use
 $F_m \assign \mbox{span} \{ E_k (x) = e^{ikx / L_x}, | k |
 \leqslant m \}$,$P_n\assign
 \mbox{span} \{ \varphi_k : 0 \leqslant k \leqslant n \}$
 as the basis set for $x$ direction and $y$ direction,
 correspondingly, where $\varphi_0 (y) = 1$,
 $\varphi_1 (y) = x$,
 $\varphi_k (y) = L_k (y) - L_{k - 2} (y)$, for $k\ge2$ and
 $L_k (y)$ denotes the Legendre polynomial of degree
 $k$. {Note that the basis set for $y$ direction is a
   direct extension of the nearly orthogonal Legendre bases
   introduced by Shen \cite{Shen1994}. Here we include
   constant and linear bases to treat Robin type boundary
   conditions. The corresponding stiffness matrix is still
   diagonal and the mass matrix is banded.}  The number of
 degree of freedom (DoF) in $x$ and $y$ direction are
 $n_x=2m+1$ and $n_y=n+1$.  For given $m$ and $n$, we take
 $H^N_\mu = H^N_\phi = F_m \otimes P_n$ as the approximation
 space for $\mu$ and $\phi$. For the Navier-Stokes equation,
 the velocity in $x$-component satisfies the GNBC, a Robin
 type boundary condition, while the component in $y$
 direction satisfies the Dirichlet boundary condition. The
 Robin type boundary condition is treated naturally in the
 weak form, the Dirichlet boundary condition is imposed on
 the approximation space. So we take the Galerkin
 approximation space for $\bu$ as
 $H^N_{\bu} = (F_m \otimes P_n) \times (F_m \otimes P_n^0)$,
 where $P_n^0 = \mbox{span} \{ \varphi_k, k = 2, ..,
 n\}$.
 The approximation space for pressure is {
   $H_p^N = F_m \otimes P_{n-2} \backslash (E_0 \otimes P_0
   )$ }.

\subsection{Solution procedure}

{The system \eqref{linear:1}-\eqref{linear:3} is a linear variable-coefficient system. }
The constant-coefficient terms in the system all lead to sparse
matrices that is time independent. However, the
variable-coefficient terms lead to time-dependent dense
matrices, thus explicitly building those time-dependent
dense matrices are extremely expensive (Note that, if one
use finite element methods, the corresponding matrices will
be sparse but still time-dependent). So we use a conjugate
gradient type solver with preconditioning (PCG), that does
not need explicitly building the matrix. Instead, it only
needs a subroutine to calculate the matrix-vector
product. Since the linear system is not symmetric, we use
BiCGSTAB method. {The preconditioning is also matrix-free. In each iteration, 
the preconditioning subroutine solve an approximated system corresponding to the system
\eqref{linear:1}-\eqref{linear:3} with convection terms in
\eqref{linear:1} and \eqref{linear:3} removed and variable coefficients
in \eqref{linear:2} approximated by constants.} {see \cite{Shen.YY2015, YuYang2017} for more details about the preconditioning.} Its effectiveness is
shown in the next section.  The solution procedure for the
BDF2 scheme is similar.

{ It worth to mention that, a new approach called
  Scalar Auxiliary Variable (SAV) method is recently
  introduced by Shen et
  al. \cite{shen_scalar_2017,shen_new_2017}. Its essential
  idea is quite similar to the IEQ type method but extending
  its applicability to more general cases. The differences
  between these two methods are that the bulk energy
  functional (integral) is quadratized in the SAV method,
  but the bulk energy density (integrand) is quadratized in
  the IEQ method. For some particular models like MBE model
  without slope selections where the nonlinear potential is
  not bounded from below but the total energy is, the SAV
  method can generates unconditionally energy stable schemes
  as well.  Moreover, the SAV method can get
  constant-coefficient terms for the gradient flow
  part. However, when considering the hydrodynamics coupled
  model, the variable-coefficient terms in the Navier-Stokes
  equation and the convection part of the phase-field
  equation are still inevitable using both methods.  }
 
\subsection{The startup step} Both CN and BDF2 schemes need
two initial steps to startup. We can use any first order
scheme to calculate $\phi^1, \bu^1, p^1, U^1, W^1$. In our
simulations, We use the first order scheme developed in
\cite{Shen.YY2015}.

\section{Numerical results}
In this section, we present various numerical experiments to
validate the developed CN scheme \eqref{cn:1}-\eqref{cnbd:6}
and BDF2 scheme \eqref{bdf:1}-\eqref{bdfbd:6}, and
demonstrate their stability and accuracy.

We examine the accuracy, stability and efficiency of the
proposed schemes by performing a classical shear flow
experiment between two parallel plates which move in
opposite directions at a constant speed. If not explicit
specified, the model parameters take default values given
below, which is consistent to the benchmark simulation in
\cite{QianWS2003,QianWS2006,GaoW2012, GaoW2014,
  Shen.YY2015,YuYang2017}.
 \begin{eqnarray}\label{para:set}
    \lambda = 20, \quad M = 0.0125, \quad \gamma = 100, 
    \quad \ell (\phi) = 1 / 0.19, \quad \nu = 1/0.6, 
    \quad  \veps = 0.05.
 \end{eqnarray}

\subsection{Convergence test for space and time}

We first test the convergence in space and time by
presenting numerical results for two cases. In case 1, we
set $\bu_w=(\pm 0.7,0)$, $\theta_s=64^\circ$, where $\bu_w$
is the velocities of upper and bottom plates. the sign of
``$\pm$'' means the values on top and bottom boundaries have
different directions, i.e., the top plate ($y=1$) moves at
the speed $0.7$ and the bottom plate ($y=-1$) moves at the
speed $-0.7$. In case 2, we set $\bu_w=(\pm 0.2,0)$ and
$\theta_s=77.6^\circ$. In both cases, $L_x=10$ and the
initial velocity field takes the profile of Couette flow,
while the initial value of $\phi$ is given as
\begin{eqnarray}\label{eq:initphi}
 \phi_0 (x, y) = \tanh \Big( \frac{1}{ \sqrt{2}\veps }\big(0.25L_x - | x - 0.5L_x | \big)\Big) . 
\end{eqnarray}
 
\begin{figure}[!htb]
  \centering
  \resizebox{0.8\columnwidth}{!}{\includegraphics{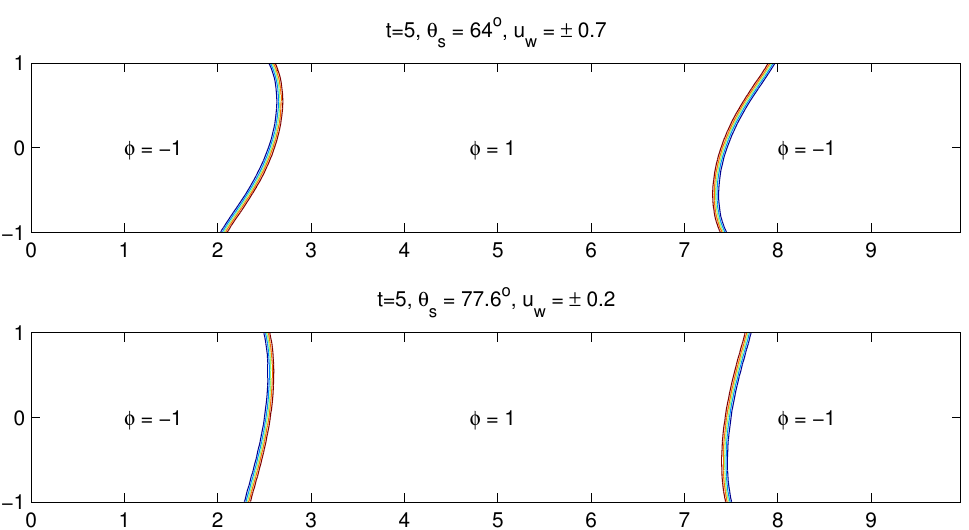}}
  \caption{\label{fig:sheartheta} The contours of phase
    variable $\phi$ at $t = 5$ that are obtained by the BDF2
    scheme with 257 Fourier modes, 32 Legendre modes and
    $\delta t=0.01$. (Top) The contour of $\phi$ for the
    case $\bu_w=(\pm 0.7, 0)$, $\theta_s=64^\circ$; (Bottom)
    The contour of $\phi$ for the case $\bu_w=(\pm 0.2, 0)$,
    $\theta_s=77.6^\circ$.  }
\end{figure}

\begin{figure}[!htb]
  \centering \subfigure[BDF2
  scheme.]{\includegraphics[width=0.48\textwidth,clip==]{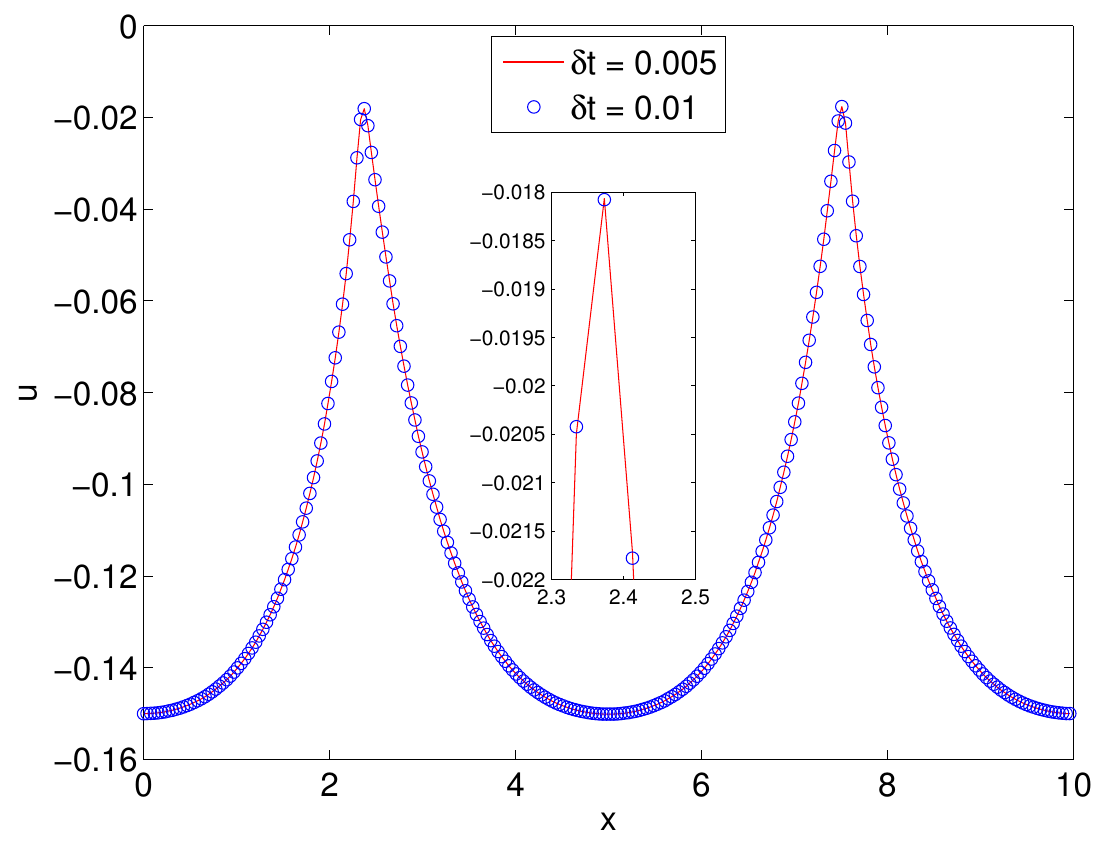}}
  \subfigure[CN
  scheme.]{\includegraphics[width=0.48\textwidth,clip==]{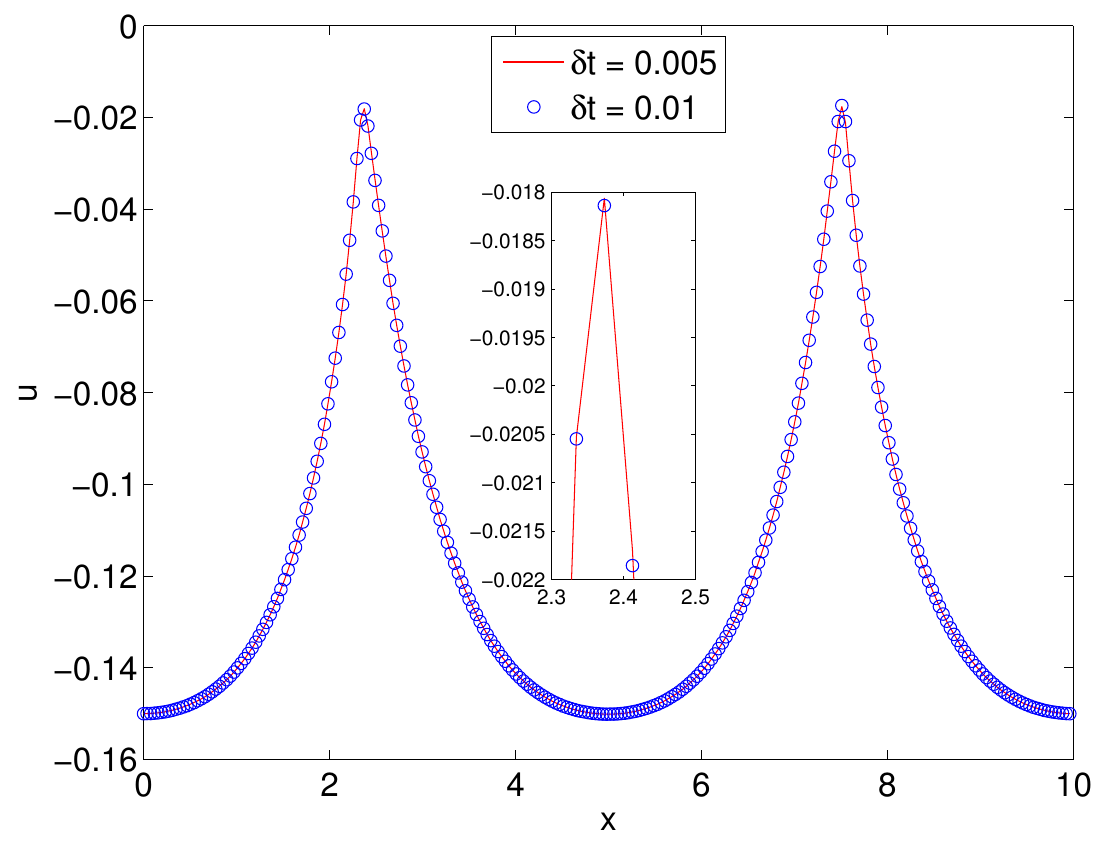}}
\caption{\label{fig:shearVx}The $x$-component velocity at
  $y=-1, t=10$ using (a) the BDF2 scheme and (b) the CN
  scheme with $\delta t=0.005$ and $\delta t=0.01$ for the
  case $\theta_s=77.6^\circ, \bu_w=(\pm 0.2,0)$.  }
\end{figure}

In Fig.\,\ref{fig:sheartheta}, we plot the contours of the
phase variable $\phi$ at $t = 5$ for the two cases, that are
obtained using the BDF2 scheme with $n_x=257$, $n_y=32$, and
$\delta t=0.01$. We only show the results of the BDF2 scheme
since the CN scheme gives visually identical results.
In Fig.\,\ref{fig:shearVx}, we plot the $x$-component of
velocity at lower boundary $y=-1$ for two different time
steps $\delta t=0.005$ and $\delta t =0.01$ obtained from
the BDF2 scheme and CN scheme for the case 2 at $t=10$. We
observe that the results obtained by two schemes are almost
identical, which means the time step $\delta t =0.01$ is
small enough to provide very accurate results for this test
case.

\begin{figure}[!htb]\centering
  \subfigure[Convergence w.r.t. the DoFs in
  $x$-direction.]{\includegraphics[width=0.49\textwidth,clip==]{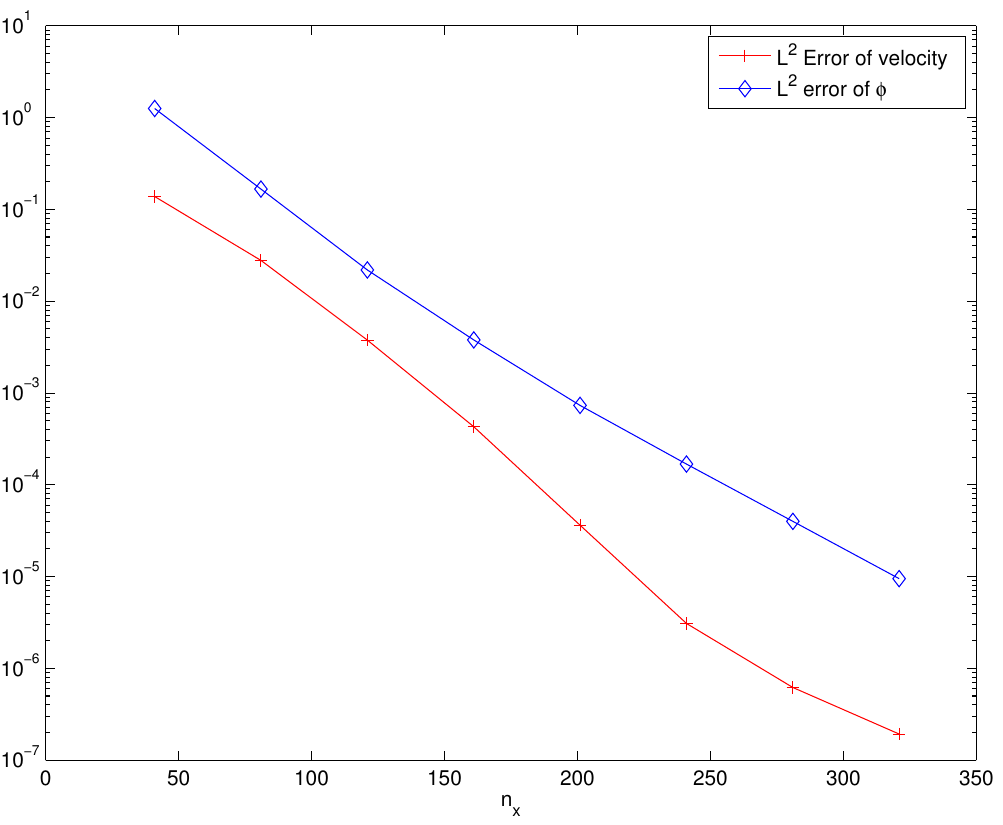}}
  \subfigure[Convergence w.r.t. the DoFs in
  $y$-direction.]{\includegraphics[width=0.49\textwidth,clip==]{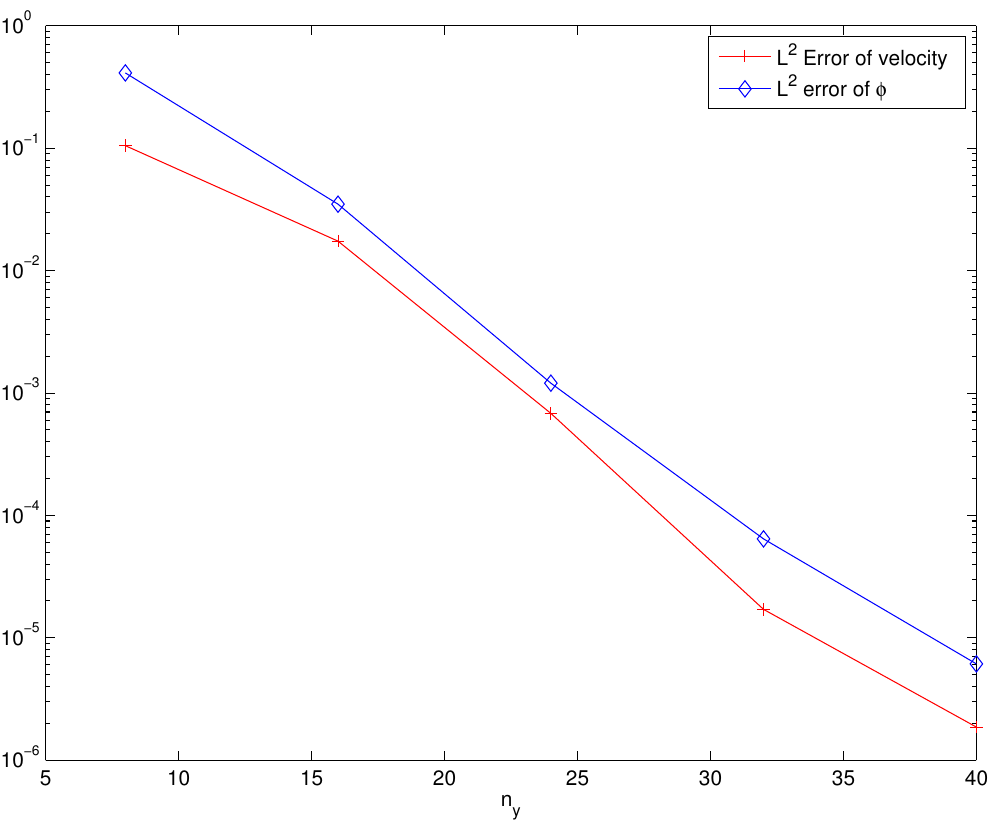}}
  \caption{\label{fig:convxy}The spatial convergence tests
    for $L^2$ error of the velocity and phase field variable
    at $T=1$ using the BDF2 scheme with the time step
    $\delta t = 0.0005$. (a): the $L^2$ error with respect
    to number of DoFs in $x$-direction; (b): the $L^2$ error
    with respect to the number of DoFs in $y$-direction.}
\end{figure}

To test the convergence for spatial discretization, we use a
very small time step $\delta t = 0.0005$ so that the errors
from the temporal discretization are negligible compared
with the spatial discretization errors.
Fig.\,\ref{fig:convxy}\,(a) shows the convergence in
$x$-direction, where we fix the number of Legendre modes
$n_y = 64$ and vary the number of Fourier modes $n_x$
starting from $41$ with the increment $40$.  The $L^2$
errors of the velocity field $\bu$ and phase variable $\phi$
are calculated at time $T = 1$, with a reference solution
obtained using the finest resolution of $n_x =511 , n_y =
64$. Similarly, Fig.\,\ref{fig:convxy}\,(b) shows the
convergence in $y$-direction, where we fix the number of
Fourier modes $n_x = 511$ for a series of $n_y$ starting
from $8$ with an increment $8$. The $L^2$ errors of the
velocity field $\bu$ and phase variable $\phi$ are again
calculated at $T = 1$ with the reference solution using the
finest resolutions of $n_x =511$ and $n_y = 64$. We see that
the proposed numerical schemes can achieve spectral accuracy
in $L^2$ norm for both velocity and phase variable.
      
\begin{figure}[h]\centering
  \subfigure[$L^2$ convergence of velocity
  $\bu$.]{\includegraphics[width=0.49\textwidth,clip==]{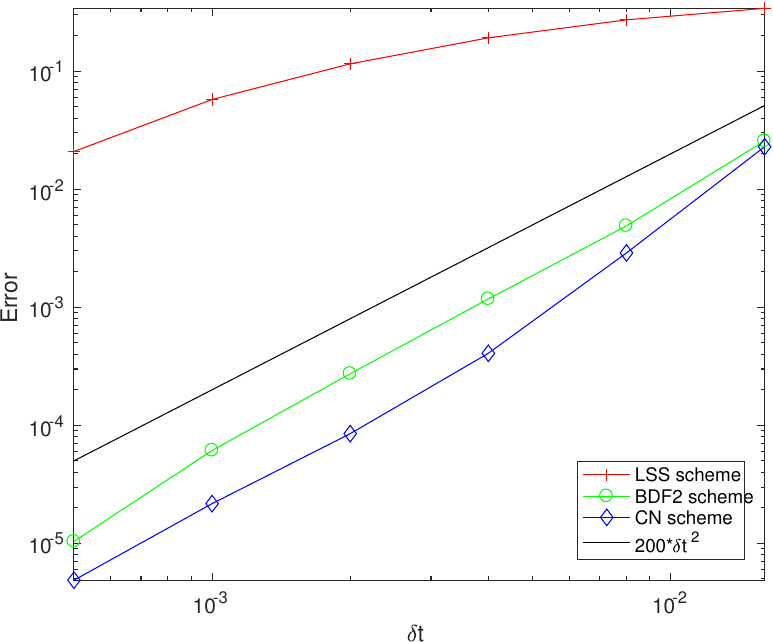}}
  \subfigure[$L^2$ convergence of phase variable
  $\phi$.]{\includegraphics[width=0.49\textwidth,clip==]{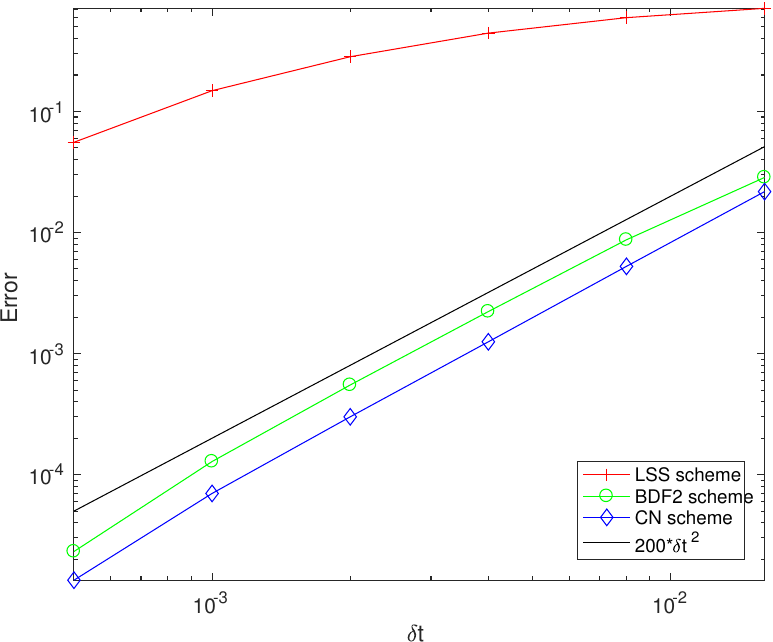}}
  \caption{\label{fig:convt}The temporal mesh refinement
    tests for the velocity field $\bu$ (a) and the phase
    field variable $\phi$ (b) obtained by CN scheme, BDF2
    scheme, and the first order scheme (denoted by LSS)
    proposed in \cite{Shen.YY2015}. The axes are in loglog
    scales, a line with slope 2 in black color is plotted to
    show the second order convergence of the CN and BDF2
    schemes.}
\end{figure}

To test the convergence for temporal discretization, we fix
the spatial resolution as $n_x= 511, n_y = 64$ so that the
errors from the spatial discretization are negligible
compared to the temporal discretization errors, and perform
the refinement test of the time step for the schemes CN and
BDF2. We choose the approximate solution using the scheme
BDF2 with time step size $\delta t=\num{2.5e-4}$ as the
benchmark solution (approximately the exact solution) for
computing errors. In Fig.\,\ref{fig:convt}, we present the
$L^2$ errors of the velocity field $\bu$ and the phase field
variable $\phi$ between the numerical solution and the exact
solution at $T= 0.4$ with different time step sizes $\delta
t = {{0.016}}/{2^{k}}$, $k = 0, 1, \ldots, 5$. The results
obtained by CN and BDF2 schemes are shown together with the
results obtained by the first order linear stabilization
scheme (denoted by LSS) proposed in \cite{Shen.YY2015} for
comparisons.  We observe that both CN and BDF2 schemes are
second order accurate for the velocity field $\bu$ as well
as the phase field variable $\phi$, which provide much more
accurate results than that of the first order LSS scheme.
From Fig.\,\ref{fig:convt}, we also observe that the accuracy of CN scheme
is better than the accuracy of BDF2 scheme, which is reasonable since
CN scheme has a smaller truncation error.
 
\subsection{Energy dissipation and volume preservation}

\begin{figure}[!htb]\centering
  \resizebox{0.49\columnwidth}{!}{\includegraphics{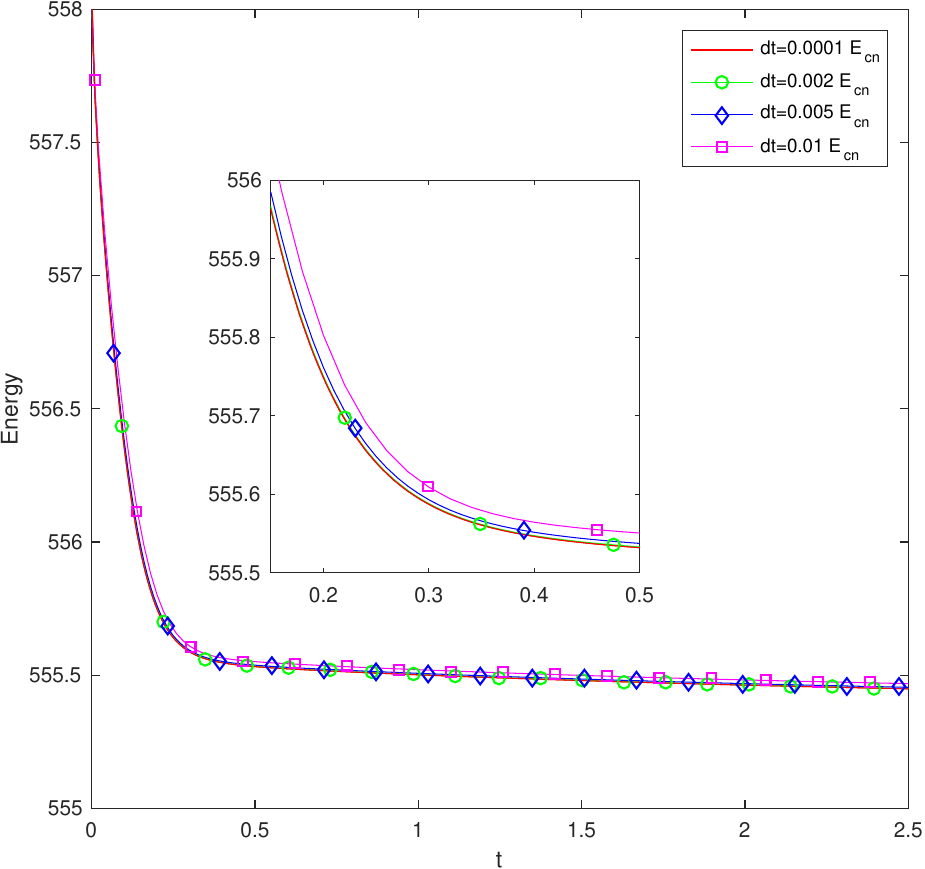}}
  \resizebox{0.49\columnwidth}{!}{\includegraphics{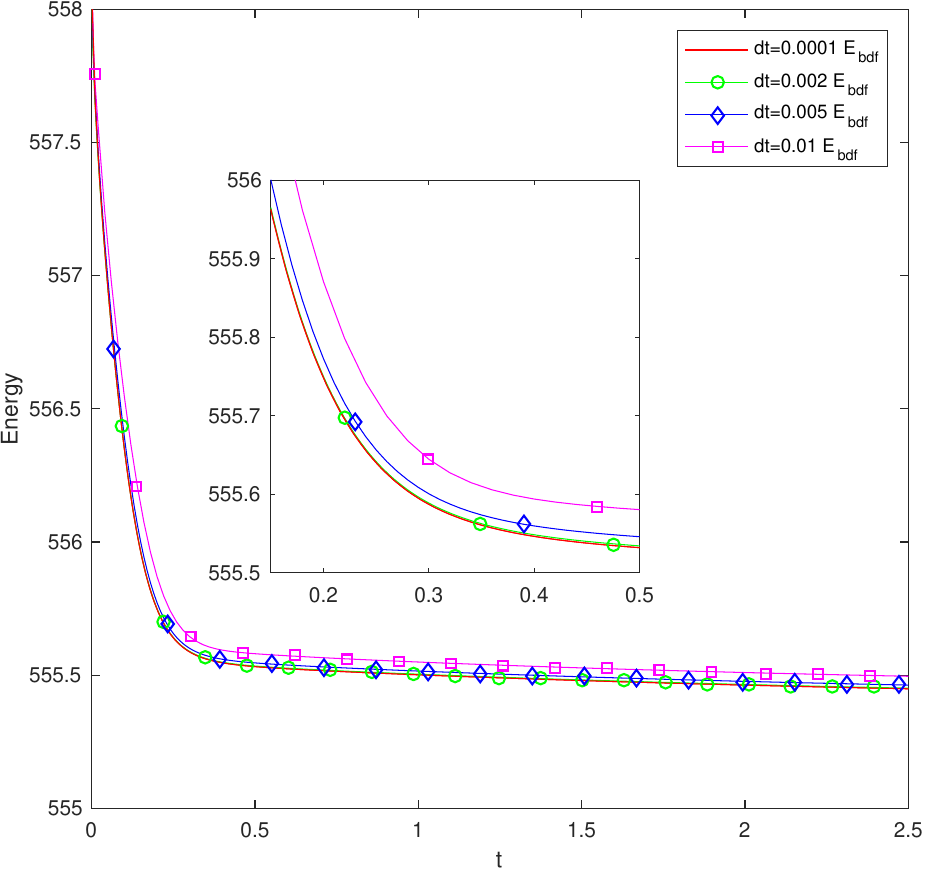}}
  \caption{\label{fig:edis} The discrete energy dissipation
    of CN (left) and BDF2 (right) schemes for several time
    steps. $u_w = 0$ and other parameters are given in
    \eqref{para:set}.}
\end{figure}

For both CN and BDF2 schemes, we test the energy
dissipation for the isolated system by setting the wall
velocity $\bu_w = 0$, and further compare the evolution of
the free energy functional for four different time step
sizes where $\delta t=0.01, 0.005, 0.002$ and $0.0001$ for
the default parameters \eqref{para:set} until $T=2.5$ in
Fig.\,\ref{fig:edis}.  For either scheme, we observe that
all energy curves show the decays, that confirm that our
algorithms are unconditionally stable. For smaller time
steps of $\delta t = 0.0001, 0.002, 0.005$, the three energy
curves coincide very well. But for the larger time step of
$\delta t=0.01$, the energy curve is considerable (but not
very far) away from others. This means the time step size
has to be smaller than $0.01$ in order to get reasonably
good accuracy.

{ To investigate the error introduced by the
  numerical approximation of \eqref{sys:5}, we plot the IEQ
  energy and the original energy for the BDF2 scheme using
  two larger time step sizes
  $\delta t= 0.01$ and $0.02$ in the left part of
  Fig. \ref{fig:ediff}. We see that the original energy still have a very good
  dissipation property, although the differences between the
  original energy and the IEQ energy are not very small. To further verify
  the convergence of the difference between the original
  energy and the modified IEQ energy, we plot in the second
  part of Fig. \ref{fig:ediff} the quantity
  $\|U \|^2 - \|\phi^2-1\|^2$ (which is the major part of
  the difference between the original and IEQ energy) at
  $t=0.4$ for both BDF2 scheme and CN scheme using different
  time steps. A clearly second order convergence is
  observed.  }

\begin{figure}[!htb]\centering
	\resizebox{0.49\columnwidth}{!}{\includegraphics{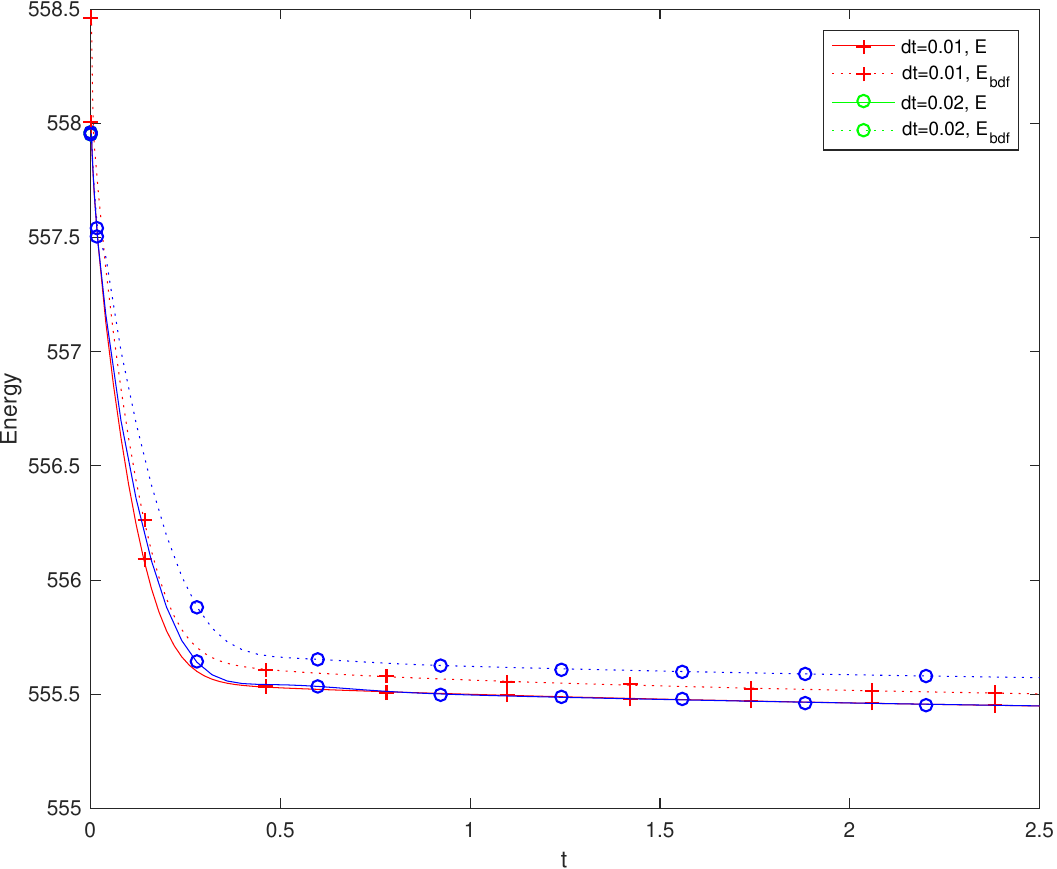}}
	\resizebox{0.49\columnwidth}{!}{\includegraphics{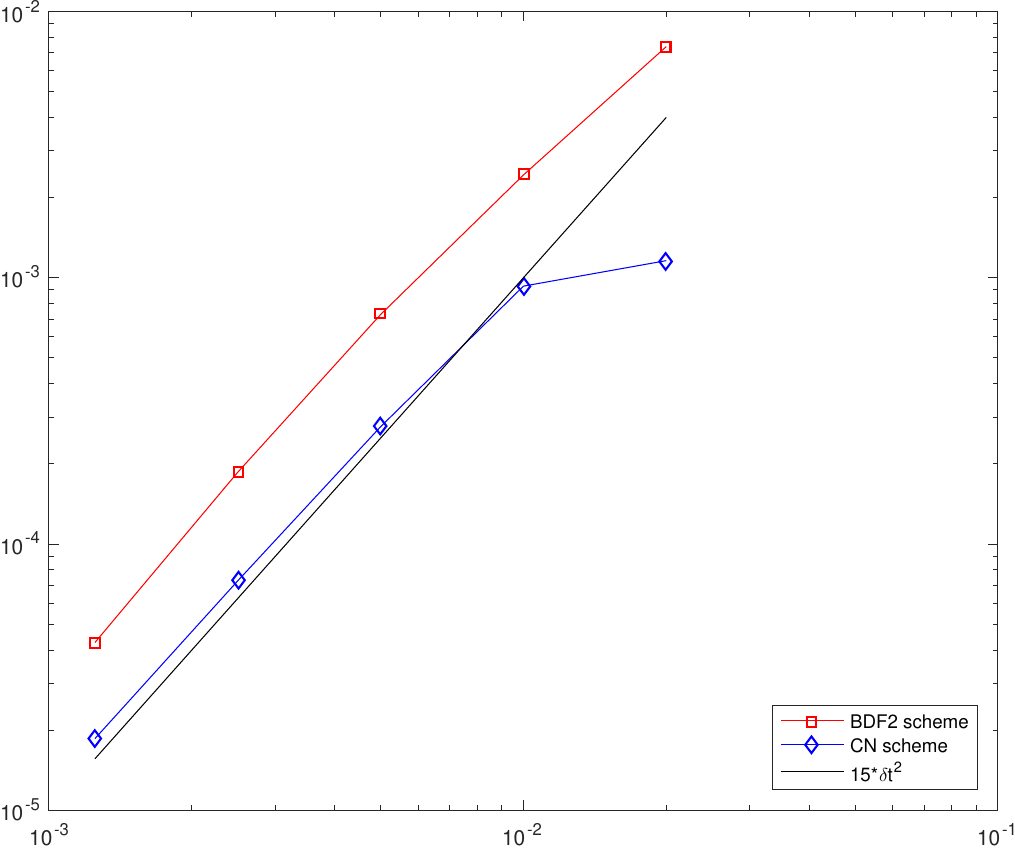}}
	\caption{\label{fig:ediff} The difference between the
      discrete BDF2 energy given by \eqref{bdf:ene:def} and the original energy given by \eqref{eq:energies} using
      different time steps. Left: The BDF2 energy (dotted lines) and
      original energy (solid lines) plot from $t=0$ to $t=2.5$ for
      different time step sizes; Right: the quantity
      $\| U \|^2 - \|\phi^2-1\|^2 $ at $t=0.4$ for BDF2
      scheme and CN scheme using different step sizes
      plotted with line $15\delta t^2$.}
	\label{fig:ediff}
\end{figure}

\begin{figure}[h]\centering
\resizebox{0.49\columnwidth}{!}{\includegraphics{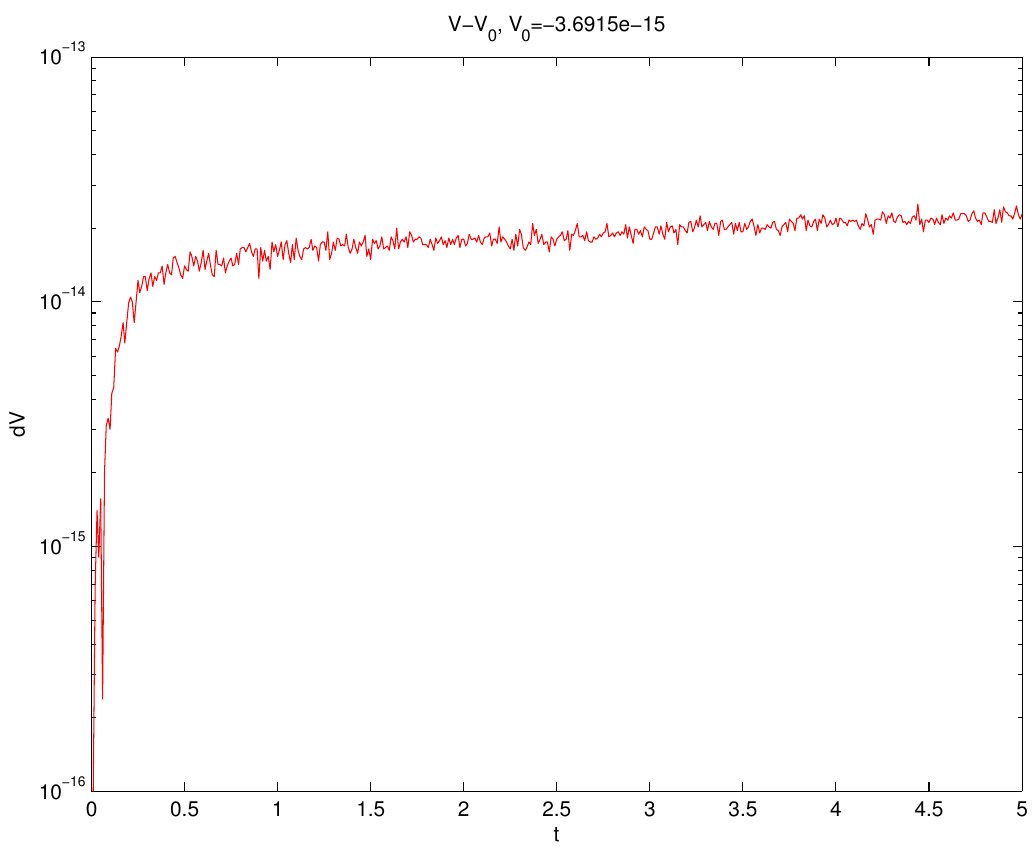}}
\resizebox{0.49\columnwidth}{!}{\includegraphics{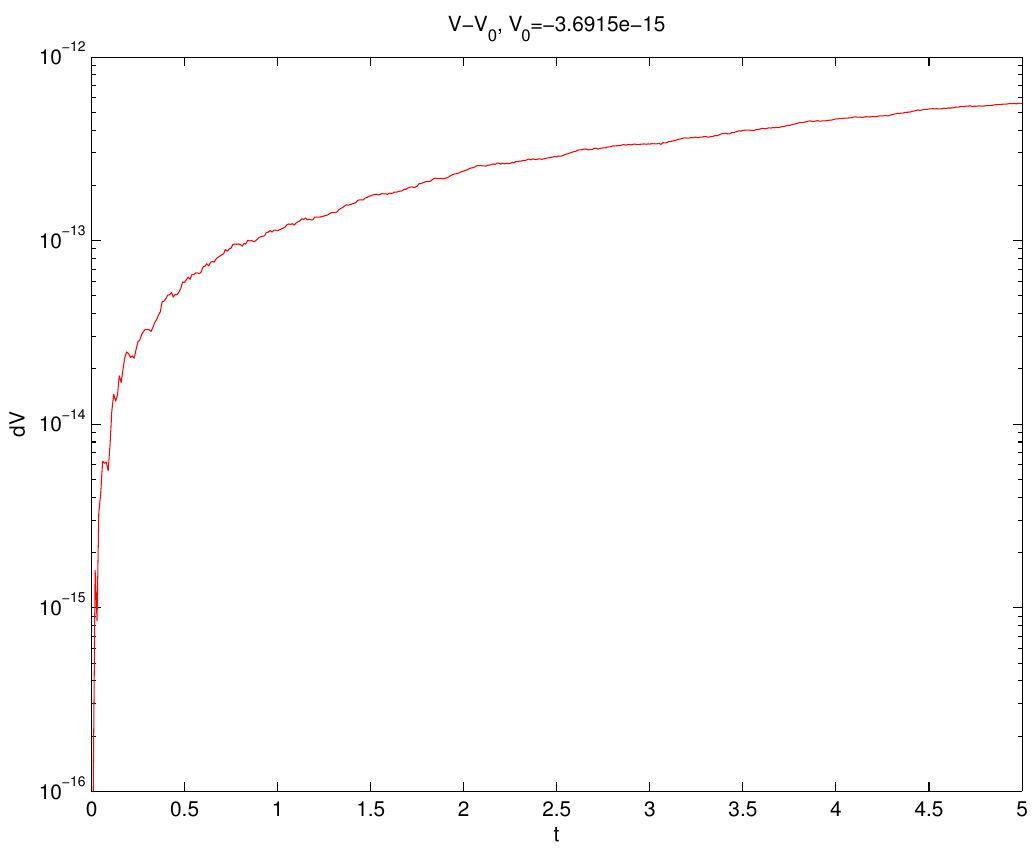}}
   \caption{Volume conservation property of CN (left) and
BDF2 (right) schemes.}\label{fig:vcons}
\end{figure}

In Fig.\ref{fig:vcons}, we show the time evolution for the
volume difference $V(t)-V_0$ where $V(t)=\int_\Omega\phi(x,
t)d\bx$ and $V_0=\int_\Omega\phi_0(x)d\bx$ for the numerical
solutions obtained by the schemes CN and BDF2, we observe
that the volume difference is very close to the machine
precision.

\subsection{Efficiency}

In Table\,\ref{tbl:n:gamma} and Table\,\ref{tbl:dt:B}, we
show the number of iterations needed by the BiCGSTAB solver
for the CN and BDF2 schemes.  The default parameters are
$n_x = 257, n_y =32$, $\delta t = 0.01$, $\gamma = 500$,
$\lambda = 12$. We vary these parameters one by one while
fixing the rest of them to be default values. In
Table\,\ref{tbl:n:gamma}, the number of iterations are
always around $O(10)$ for various grid points and $\gamma$,
that means the number of iterations actually does not show
any dependence on the number of spatial grid points and the
parameter $\gamma$. In Table\, \ref{tbl:dt:B}, when we vary
the time step $\delta t$ and parameter $\lambda$, we can see
that larger values of them can cause significant increase
for the number of iterations. Namely, when $\delta t $ and
$\lambda$ are larger, the conditional numbers of the linear
system become worse that is reasonable and can be easily
observed from the form of the linear system
\eqref{rewrite:cn1}.
\begin{table}[!htb]\centering
     \begin{tabular}{p{9ex}<{\centering} | *{3}{p{7ex}<{\centering}} }
      \hline
      $n_x \times n_y$ & 129$\times$16 & 257$\times$32 & 513$\times$64\\
      \hline
      CN & 8 & 8.7 & 9.6\\
      \hline
      BDF2 & 8.5 & 8.5 & 8.5\\
      \hline
     \end{tabular}\quad 
	\begin{tabular}{p{9ex}<{\centering} | *{3}{p{6ex}<{\centering}} }
	  \hline
	  $\gamma$ & 100 & 10 & 1\\
	  \hline
	  CN & 7.9 & 9 & 9\\
	  \hline
	  BDF2 & 9.1 & 10.0 & 10.1\\
	  \hline
	\end{tabular}\vskip 0.5cm
    \caption{\label{tbl:n:gamma} 
      The average number of the inner iterations for BiCGSTAB with respect to the grid points 
       $n_x\times n_y$ and the parameter $\gamma$. The tolerance is  
       $10^{-8}$ for schemes CN and BDF2. The default values are 
         $n_x = 257, n_y = 32$, $\delta t = 0.01$, $\gamma = 500$, $\lambda = 12$. }
\end{table}

\begin{table}[!htb]\centering
     \begin{tabular}{p{9ex}<{\centering} | *{3}{p{7ex}<{\centering}} }
     	\hline
       $\mbox{$\delta t$}$ & 0.001 & 0.1 & 1\\
       \hline
       CN & 3.8 & 27 & 68\\
       \hline
       BDF2 & 4 & 27 & 69.6\\
       \hline
     \end{tabular}\quad
    \begin{tabular}{p{9ex}<{\centering} | *{3}{p{6ex}<{\centering}}}
       \hline
       $\lambda$ & 1 & 60 & 144\\
       \hline
       CN & 5 & 19.6 & 35.5\\
       \hline
       BDF2 & 5 & 20.8 & 42.5\\
       \hline
     \end{tabular}\vskip 0.5cm
     \caption{\label{tbl:dt:B} 
       The average number of the inner iterations for BiCGSTAB  with respect to the time step 
       $\delta t$ and the parameter $\lambda$. The tolerance  
       is $10^{-8}$ for both schemes CN and BDF2. The default values are 
       $n_x = 257, n_y =32$, $\delta t = 0.01$, $\gamma = 500$, $\lambda = 12$. }
\end{table}

\subsection{Simulation of a drop in shear flow}
In this subsection, we simulate the dynamics of a drop in
shear flow using BDF2 scheme. The channel length
$L_x=10$. The wall velocity $\bu_w=(\pm2,0)$. All other
parameters are given in \eqref{para:set}.

\begin{figure}[!htb]\centering
  \resizebox{0.48\columnwidth}{!}{\includegraphics{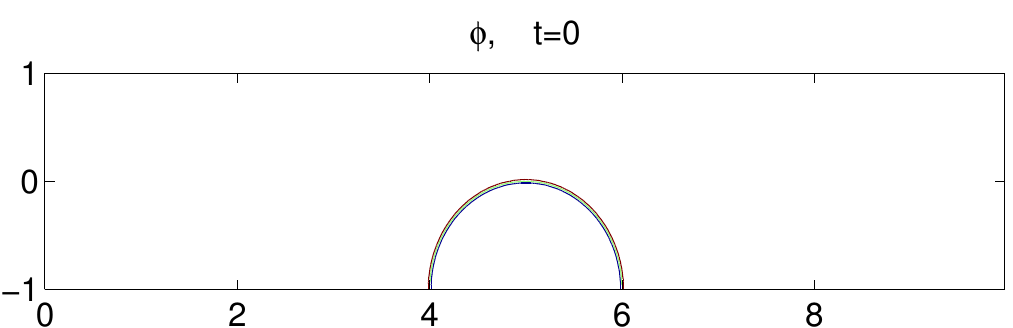}}
  \resizebox{0.48\columnwidth}{!}{\includegraphics{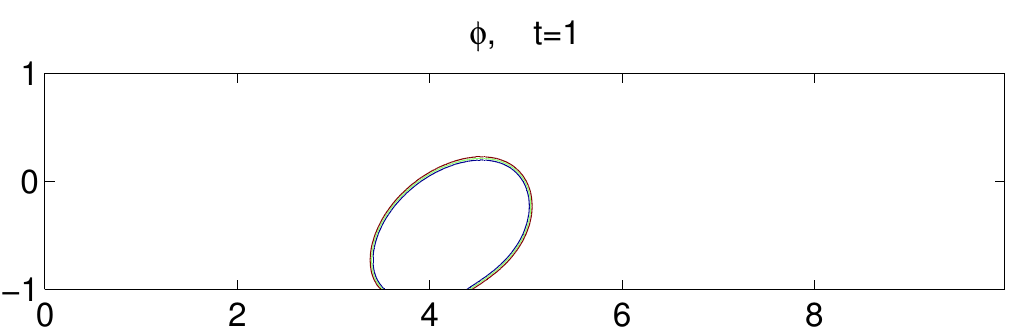}}\\
  \resizebox{0.48\columnwidth}{!}{\includegraphics{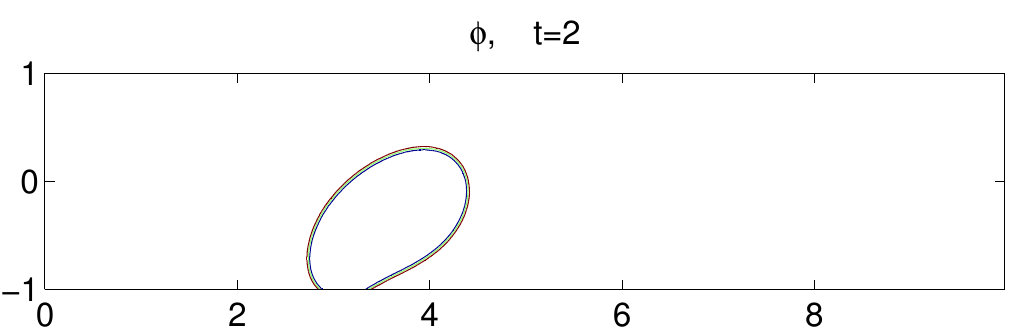}}
  \resizebox{0.48\columnwidth}{!}{\includegraphics{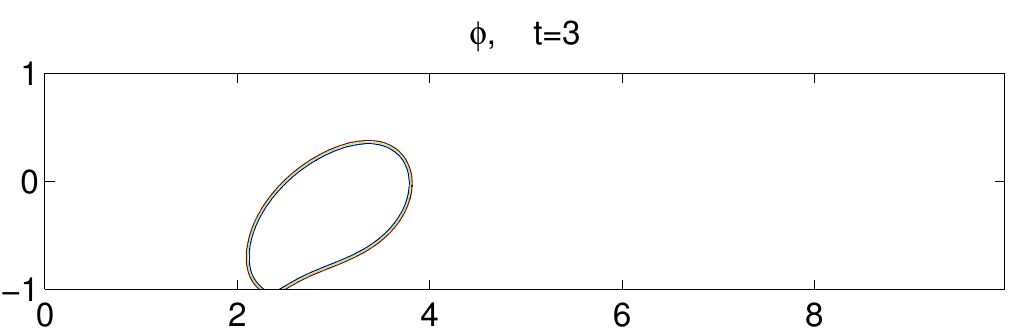}}\\
  \resizebox{0.48\columnwidth}{!}{\includegraphics{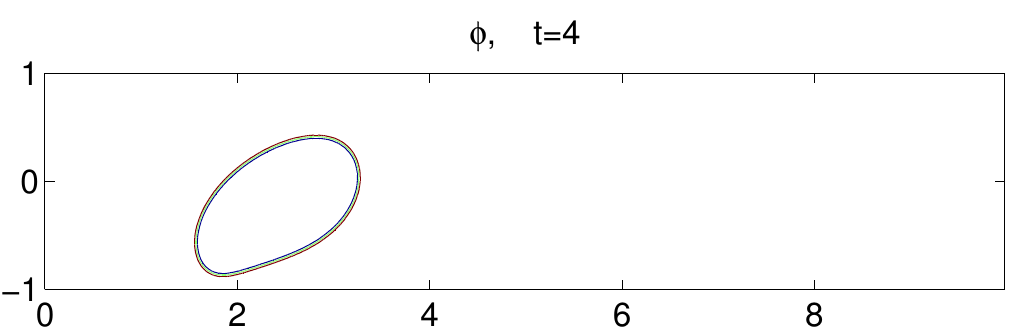}}
  \resizebox{0.48\columnwidth}{!}{\includegraphics{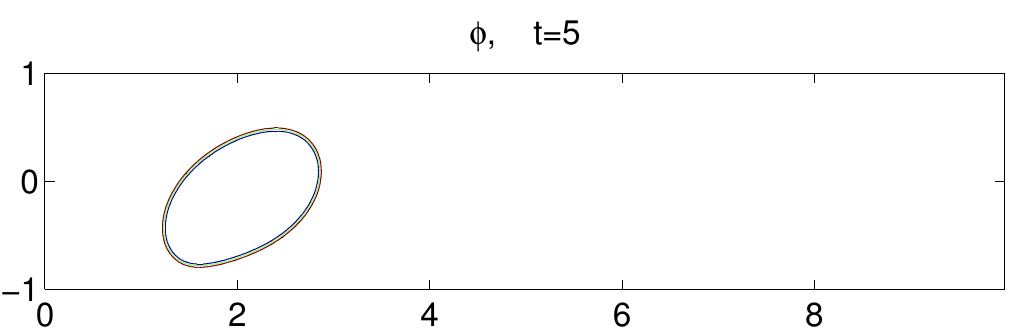}}
  \caption{The dynamical behaviors of a drop in shear flow
    simulated using the BDF2 scheme. Snapshots are taken at
    $t=0,1,2,3,4,5$. We take the wall velocity $\bu_w=(\pm
    2,0)$ and the static contact angle $\theta_s =
    {\pi}/{6}$.}\label{fig:dropshear}
\end{figure}

\begin{figure}[!htb]\centering
  \resizebox{0.48\columnwidth}{!}{\includegraphics{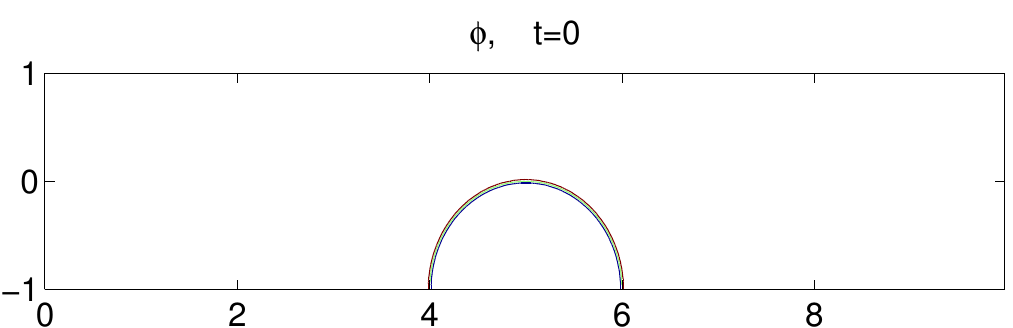}}
  \resizebox{0.48\columnwidth}{!}{\includegraphics{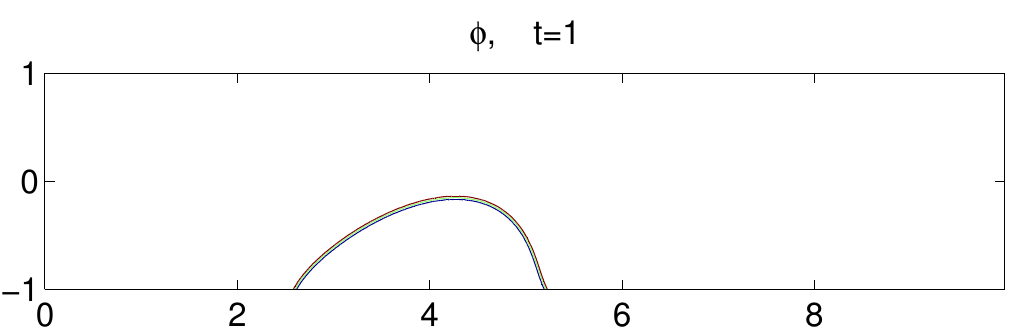}}\\
  \resizebox{0.48\columnwidth}{!}{\includegraphics{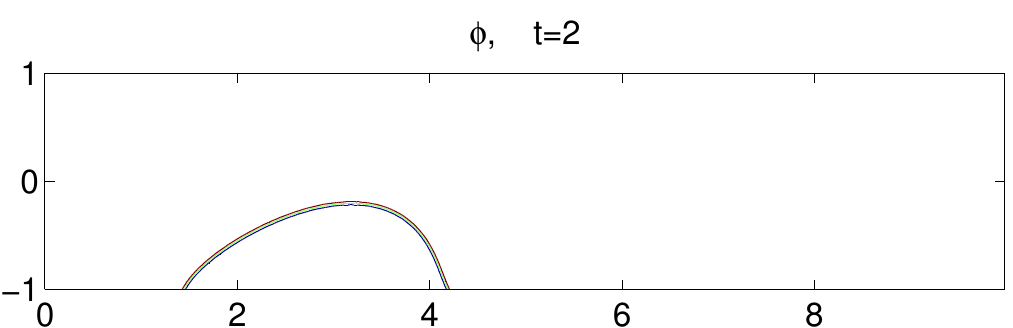}}
  \resizebox{0.48\columnwidth}{!}{\includegraphics{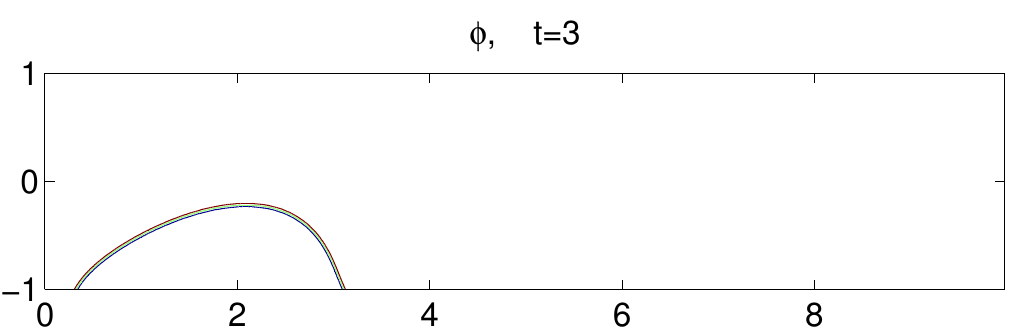}}\\
  \resizebox{0.48\columnwidth}{!}{\includegraphics{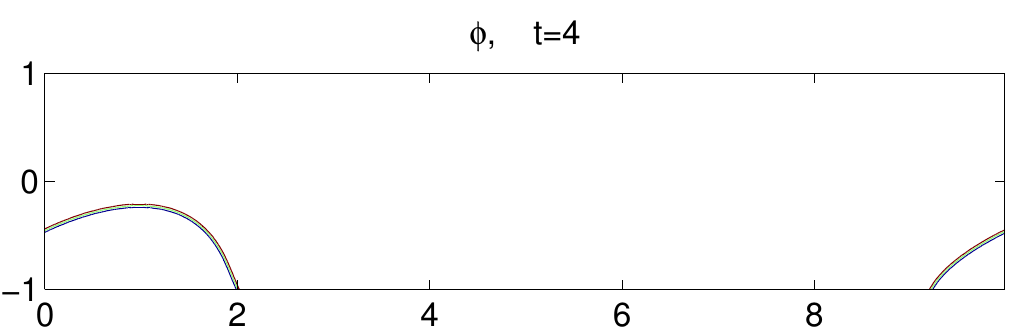}}
  \resizebox{0.48\columnwidth}{!}{\includegraphics{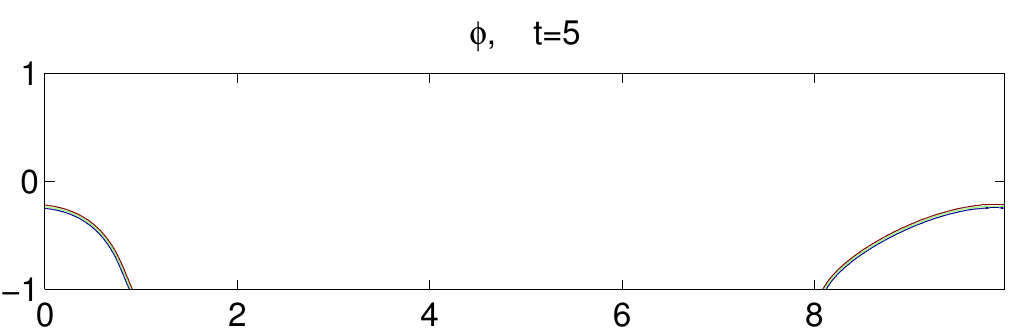}}
  \caption{The dynamical behaviors of a drop in shear flow
    using the BDF2 scheme. Snapshots are taken at
    $t=0,1,2,3,4,5$. We take the wall velocity $\bu_w=(\pm
    2,0)$ and the static contact angle $\theta_s =
    {2\pi}/{3}$.}\label{fig:dropshear2}
\end{figure}

Fig.\,\ref{fig:dropshear} illustrates the dynamical motions
of a drop under shear with an acute contact angle
$\theta_s={\pi}/{6}$ until $T=5$.  Initially, the drop is
set in the middle of the bottom plate, as shown in the first
subfigure. As the bottom plate moves, the drop also moves
with the plate but at a much smaller velocity. This means
the drop slips on the bottom plate. As time goes on, the
contact region of the drop with the bottom boundary gets
smaller and smaller. It eventually gets off the bottom
boundary around $t=3$ and moves toward the center of the
channel.  We also simulate an obtuse contact angle $\theta_s
= {2\pi}/{3}$ case, in which the drop is harder to get off
the bottom boundary. The result is shown in
Fig.\,\ref{fig:dropshear2}.

 \section{Concluding remarks}

 By combining the projection method for Navier-Stokes
 equations, the IEQ method for the nonlinear bulk and
 boundary energy gradients, and a subtle explicit-implicit
 technique for the stress and convective terms, we have
 constructed two linear, second order unconditionally stable
 temporal discretization schemes for the phase field MCL
 model.  The well-posedness of the semi-discretized linear
 systems and their energy stabilities are proved rigorously.
 Numerical simulations have also verified both schemes are
 unconditionally stable and second order accurate, while the
 CN scheme behaves a little bit better than the BDF2 scheme.
 Although we have considered only time discretizations here,
 the results can carry over to any consistent
 finite-dimensional Galerkin (finite element or spectral)
 approximations since the proofs are all based on
 variational formulations with all test functions in the
 same space as the trial function.

 \section*{Acknowledgments}
 
 The work of X. Yang was partially supported by the National
 Science Foundation under grant number NSF
 DMS-1418898 and DMS-1720212. The work of H. Yu was partially supported by
 NSFC projects 11771439, 91530322 and China National Program on Key Basic Research Project 2015CB856003.

\bibliographystyle{plain}
%\bibliographystyle{elsarticle-num} 
%\bibliography{references2}

\end{document}